\documentclass[12pt]{amsart}
\usepackage[T1]{fontenc}
\usepackage[cp1250]{inputenc}
\usepackage[dvips]{graphicx}
\usepackage{amssymb}
\usepackage{amsmath}
\usepackage{amsthm}
\usepackage{amstext}
\usepackage{amsopn}
\usepackage{mathrsfs}
\usepackage{latexsym}
\usepackage{mathabx}
\allowdisplaybreaks
\newtheorem{Definition}{Definition}[section]
\newtheorem{Proposition}{Proposition}[section]
\newtheorem{Lemma}{Lemma}[section]
\newtheorem{Theorem}{Theorem}[section]
\newtheorem{Corollary}{Corollary}[section]
\newtheorem{Remark}{Remark}[section]
\newtheorem{Example}{Example}[section]

\numberwithin{equation}{section}
\begin{document}
%%%%%%%%%%%%%%%%%%%%%%%%%%%%%%%%%%%%%%%%%%%%%%%%%%%%%%%%%%%%
\bibliographystyle{plain}
\footnotetext{
This work is partially supported by MNiSW research grant No 1 P03A 013 30}
%%%%%%%%%%%%%%%%%%%%%%%%%%%%%%%%%%%%%%%%%%%%%%%%%%%%%%%%%%%%%%%%%
\title[Matricially free random variables]
{Matricially free random variables}
\author[R. Lenczewski]{Romuald Lenczewski}
\address{Romuald Lenczewski, \newline
Instytut Matematyki i Informatyki, Politechnika Wroc\l{}awska, \newline
Wybrze\.{z}e Wyspia\'{n}skiego 27, 50-370 Wroc{\l}aw, Poland  \vspace{10pt}}
\email{Romuald.Lenczewski@pwr.wroc.pl}
%%%%%%%%%%%%%%%%%%%%%%%%%%%%%%%%%%%%%%%%%%%%%%%%%%%%%%%%%%%%%%%%
\maketitle
\begin{abstract}
We show that the operatorial framework developed by Voiculescu for free random variables
can be extended to arrays of random variables whose multiplication
imitates matricial multiplication. The associated notion of independence, called {\it matricial freeness}, 
can be viewed as a generalization of both freeness and monotone independence. 
At the same time, the sums of matricially free random variables, called {\it random pseudomatrices}, 
are closely related to Gaussian random matrices. 
The main results presented in this paper concern the standard and tracial central 
limit theorems for random pseudomatrices and the corresponding limit distributions which
can be viewed as matricial generalizations of semicirle laws.
\end{abstract}

\section{Introduction}
It has been shown by Voiculescu [22] that free random variables arise naturally as limits of random matrices. In particular, if we take symmetric matrices whose entries form a family of independent Gaussian random variables and we let the size of these matrices go to infinity, their moments (with respect to normalized trace composed with classical expectation) converge to the moments of freely independent random variables with the semicirle distribution obtained by Wigner [25] as the limit distribution of one Gaussian random matrix.

Therefore, we can study free random variables, using at least two different frameworks: 
operator algebras and random matrices. However, it is to some extent surprising that 
a random matrix framework, of quite different nature than that of operator algebras, 
exists for free random variables, and the connection between these two approaches does not seem to be very transparent. In this connection, our first motivation is to better understand the relation between the operatorial approach to free probability and random matrices. 

The second motivation comes from the question whether different types of independence, like freeness and monotone independence
of Muraki [16], can be included in one natural model. Note in this context that models with more than one state on a given algebra, like conditional freeness of Bo\.{z}ejko and Speicher [6] and freeness with infinitely many states of Cabanal-Duivillard and Ionesco [9,10] extend free probability and include, as shown by Franz [11], certain elements of monotone probability. For instance, this can be done for convolutions, but including monotone independence in the framework of conditional freeness can be done only under additional (rather restrictive) assumptions on the considered algebras. We show in this paper that one can remedy this situation by introducing a concept of `independence' which reminds freeness, but at the same time has some `matricial' features which places it somewhere between freeness and the model of random matrices.

A different reason to look for a new concept of independence arises from concrete examples of interpolations between free probability and monotone probability [14,15]. In particular, the continuous $(p,q)$-Brownian motions with Kesten distributions and related Poisson processes studied in [2,15] lead to the first example of a two-mode interacting Fock space introduced directly and not by means of orthogonal polynomials. This example has certain `matricial' features which also call for a new model of `independence' that would be related to freeness.

The main result of our paper is the construction of a model,
called {\it matricial freeness}, which is related to the concept of 
the free product of states introduced and studied by Ching [10] in the context of von Neumann algebras 
and by Avitzour [3] and Voiculescu [20] in the context of $C^{*}$-algebras. 
The underlying concept is that of the matricially free product of an array of Hilbert spaces with distinguished unit 
vectors $({\mathcal H}_{i,j}, \xi_{i,j})$ which reminds the free product of Hilbert spaces.
However, it also satisfies the condition imitating matrix multiplication 
and the condition of `diagonal subordination' which says that tensor products must end with
`diagonal' Hilbert spaces. Similarities between free probability and 
`matricially free probability' hold also on other levels, some of which 
are studied in this paper. Strictly speaking, however, one needs to take a restriction of the 
matricially free product of states, called the {\it strongly matricially free product} 
of states, to recover the free product as a special case. 

Roughly speaking, if $(X_{i,j}(n))_{1\leq i,j \leq n}$ is an array of 
self-adjoint matricially free random variables in a *-algebra ${\mathcal A}_{n}$, equipped
with a distinguished state $\phi_{n}$ and a sequence $(\psi_{n,j})_{1\leq j \leq n}$ 
of additional states, called `conditions', for each natural $n$, 
we can take the sums
$$
S(n)=\sum_{i,j=1}^{n}X_{i,j}(n)
$$
called {\it random pseudomatrices}, and study their asymptotic distributions (as $n\rightarrow \infty$) with respect to 
the states $\phi_{n}$ and with respect to normalized traces
$$
\psi_{n}=\frac{1}{n}\sum_{j=1}^{n}\psi_{n,j},
$$
respectively. We assume that the distributions of the $X_{i,j}(n)$ in 
the states $\phi_n$ and $\psi_{n,k}$ are 
not identical and depend on $n$ in such a way that their 
$\psi_{n}$-distributions remind those of Gaussian random matrices
in the approach of Voiculescu. This relation to random matrices 
as well as asymptotic matricial freeness of random pseudomatrices will be studied in a forthcoming paper. 
In progress is the work on the addition of matricially free random variables and the resulting matricial R-transform.

It turns out that the limit theorem for the above random pseudomatrices 
in the states $\phi_n$, which we call `standard', may be viewed as an analog of the central limit theorem for free random 
variables (especially, if we take square arrays). In turn, the limit theorem for random pseudomatrices in 
the states $\psi_n$, called `tracial', is related to the limit theorem for random matrices 
(especially, if we take square arrays). 
The limit distributions play then the role of multivariate generalizations of the 
semicircle distributions. Let us point out, however, that when we consider
triangular arrays, the framework of matricial freeness can as well be viewed 
as a generalization of that of monotone independence.

In Section 2, we introduce the concepts of the `matricially free product of states' and
the `matricially free Fock space'. We obtain from these structures 
their strong counterparts in Section 3.
In Section 4, we introduce the notions
of `matricial freeness' and `strong matricial freeness' and discuss 
the example of the discrete (strongly) matricially free Fock space.
In Section 5, of combinatorial nature, we define and study certain 
real-valued functions on the set of non-crossing partitions, defined in terms of
traces of certain matrices.
In Section 6, we study the asymptotic behavior of random pseudomatrices 
and we prove standard and tracial central limit theorems.
The limit distributions, which can be interpreted as matricial multivariate generalizations
of semicirle laws, are studied in Section 7.
Their decompositions in terms of s-free additive convolutions in the case of two-dimensional arrays are proved in Section 8. 
Two geometric realizations of the limit distributions, in terms of 
walks on weighted binary trees and in terms of weighted Catalan paths, 
are given in Section 9.

\section{Matricially free products}
In this Section we introduce the notion of the matricially free product of states
as well as the corresponding notions of the matricially free product of Hilbert spaces 
and the matricially free Fock space.

When speaking of arrays indexed by two indices, say $i,j$, we shall usually assume 
that $i,j\in I$, where $I$ is an index set.  This refers to the situation when we deal with
square arrays. However, we also want to consider other arrays, like triangular arrays of the form
$T:=\{(k,l):k\geq l, ;k,l\in I\}$, where $I$ is a linearly ordered index set.
Therefore, by an {\it array} we will understand a subarray of a square array which 
{\it includes the diagonal}. Without loss of generality we can use the square array 
formulation most of the time. Of special interest will be the 
finite-dimensional case when $I=[n]:=\{1,2,\ldots , n\}$.

\begin{Definition}
{\rm Let $\widehat{{\mathcal H}}:=({\mathcal H}_{i,j})$ 
be an array of complex Hilbert spaces.
By the {\it matricially free Fock space over $\widehat{\mathcal H}$}
we understand the Hilbert space direct sum
$$
\mathcal{M}(\widehat{\mathcal H})=
{\mathbb C}\Omega \oplus \bigoplus _{m=1}^{\infty}
\bigoplus_{\stackrel {(i_{1},i_2)\neq \ldots \neq (i_{m},i_m)}
{\scriptscriptstyle n_1, \ldots , n_m \in {\mathbb N}}}
{\mathcal H}_{i_1,i_2}^{\otimes \,n_1}
\otimes 
{\mathcal H}_{i_2,i_3}^{\otimes \,n_2}
\otimes \ldots \otimes 
{\mathcal H}_{i_{m},i_{m}}^{\otimes \,n_m}
$$
where $\Omega $ is a unit vector, with the canonical inner product.}
\end{Definition}

Let us observe that the Hilbert space tensor powers which appear 
in the above direct sum have the following three properties:
\begin{enumerate}
\item the `matricial property' -- the second index of the 
preceding power agrees with the first index of the following power,
\item the `freeness property' -- the consecutive pairs of indices are different,
\item the `diagonal subordination property' -- the last pair
is `diagonal'.
\end{enumerate}

Of course, if the index set $I$ consists of one element and thus 
$\widehat{\mathcal H}$ is just one Hilbert space ${\mathcal H}$, the corresponding 
matricially free Fock space reduces to the usual free Fock space ${\mathcal F}({\mathcal H})$.
In general, however, $\mathcal{M}(\widehat{\mathcal H})$ is a (usually, proper) subspace of the free Fock space
$\mathcal{F}(\bigoplus_{i,j}{\mathcal H}_{i,j}) \cong *_{i,j}\mathcal{F}({\mathcal H}_{i,j})$.

Related to the `matricially free Fock space' is the `matricially free product of Hilbert
spaces'. The terminology parallels that introduced in free probability [20,23].
\begin{Definition}
{\rm Let $({\mathcal H}_{i,j},\xi_{i,j})$ be an array of Hilbert spaces
with distinguished unit vectors. By the {\it matricially free product 
of $({\mathcal H}_{i,j},\xi_{i,j})$} we understand the pair $({\mathcal H}, \xi)$, where
$$
{\mathcal H}= {\mathbb C}\xi \oplus\bigoplus_{m=1}^{\infty}
\bigoplus_{(i_{1},i_2)\neq \ldots \neq (i_{m},i_m)}
{\mathcal H}_{i_1,i_2}^{0}
\otimes 
{\mathcal H}_{i_2,i_3}^{0}
\otimes \ldots \otimes 
{\mathcal H}_{i_{m},i_{m}}^{0},
$$
with ${\mathcal H}_{i,j}^{0}={\mathcal H}_{i,j}\ominus {\mathbb C} \xi_{i,j}$ and
$\xi $ being a unit vector. We denote it 
$({\mathcal H}, \xi)=*^{M}_{i,j}({\mathcal H}_{i,j},\xi_{i,j})$.}
\end{Definition}

\begin{Proposition}
It holds that 
\begin{eqnarray*}
(\mathcal{M}(\widehat{\mathcal H}),\xi)
&\cong& 
*^{M}_{i,j}
({\mathcal M}({\mathcal H}_{i,j}),\xi_{i,j})
\end{eqnarray*}
\end{Proposition}
{\it Proof.}
We use the definition of the matricially free Fock space, the isomorphism
${\mathcal M}({\mathcal H}_{i,j})\cong {\mathcal F}({\mathcal H}_{i,j})$ for any $i,j$ and regroup terms.
\hfill $\blacksquare$\\
\indent{\par}
For any $j$, introduce {\it diagonal subspaces} of ${\mathcal H}$ of the form
$$
{\mathcal H}(j,j)={\mathbb C}\xi\oplus \bigoplus_{m=2}^{\infty}
\bigoplus_{
{\stackrel 
{(j,i_2)\neq \ldots \neq (i_{m},i_m)}
{\scriptscriptstyle i_2\neq j}
}}
{\mathcal H}_{j,i_2}^{0}
\otimes 
{\mathcal H}_{i_2,i_3}^{0}
\otimes \ldots \otimes 
{\mathcal H}_{i_{m},i_{m}}^{0},
$$
and the associated {\it diagonal partial isometries} $V_{j,j}:{\mathcal H}_{j,j}\otimes {\mathcal H}(j,j)\rightarrow {\mathcal H}$:
\begin{eqnarray*}
\xi_{j,j}\otimes \xi &\rightarrow &\xi\\
{\mathcal H}_{j,j}^{0}\otimes \xi &\rightarrow & {\mathcal H}_{j,j}^{0}\\
\xi_{j,j}\otimes ({\mathcal H}_{j,j_1}^{0}
\otimes \ldots \otimes {\mathcal H}_{j_{m},j_{m}}^{0})
&\rightarrow & 
{\mathcal H}_{j,j_1}^{0}
\otimes \ldots \otimes {\mathcal H}_{j_{m},j_{m}}^{0}\\
{\mathcal H}_{j,j}^{0}\otimes ({\mathcal H}_{j,j_1}^{0}
\otimes \ldots \otimes {\mathcal H}_{j_{m},j_{m}}^{0})
&\rightarrow & 
{\mathcal H}_{j,j}^{0}\otimes
{\mathcal H}_{j,j_1}^{0}
\otimes \ldots \otimes {\mathcal H}_{j_{m},j_{m}}^{0}
\end{eqnarray*}
where $m>1$.

For any $i\neq j$ we introduce {\it non-diagonal subspaces} of $\mathcal{H}$ of the form
$$
{\mathcal H}(i,j)=\bigoplus_{m=1}^{\infty}
\bigoplus_{(j,i_2)\neq \ldots \neq (i_{m},i_m)}
{\mathcal H}_{j,i_2}^{0}
\otimes 
{\mathcal H}_{i_2,i_3}^{0}
\otimes \ldots \otimes 
{\mathcal H}_{i_{m},i_{m}}^{0}
$$
and the associated {\it non-diagonal partial isometries} $V_{i,j}:{\mathcal H}_{i,j}\otimes {\mathcal H}(i,j)\rightarrow {\mathcal H}$ for $i\neq j$:
\begin{eqnarray*}
\xi_{i,j}\otimes ({\mathcal H}_{j,j_1}^{0}
\otimes \ldots \otimes {\mathcal H}_{j_{m},j_{m}}^{0})
&\rightarrow & 
{\mathcal H}_{j,j_1}^{0}
\otimes \ldots \otimes {\mathcal H}_{j_{m},j_{m}}^{0}\\
{\mathcal H}_{i,j}^{0}\otimes ({\mathcal H}_{j,j_1}^{0}
\otimes \ldots \otimes {\mathcal H}_{j_{m},j_{m}}^{0})
&\rightarrow & 
{\mathcal H}_{i,j}^{0}\otimes
{\mathcal H}_{j,j_1}^{0}
\otimes \ldots \otimes {\mathcal H}_{j_{m},j_{m}}^{0},
\end{eqnarray*}
where $m\geq 1$. 

Each ${\mathcal H}(i,j)$ is spanned by simple tensors which do not begin with vectors from 
${\mathcal H}_{i,j}^{0}$ and for that reason it is suitable for the left free action of the operators creating such vectors. Thus, roughly speaking, both types of partial isometries jointly
replace the unitary maps used in free probability. It is the diagonal subordination
property which is responsible for distinguishing two types of isometries.

Consider an array of $C^{*}$-algebras $({\mathcal A}_{i,j})$, each with a unit $1_{i,j}$ and a state $\varphi_{i,j}$, and let 
$({\mathcal H}_{i,j},\pi_{i,j},\xi_{i,j})$ be the associated GNS triples, so that
$\varphi_{i,j}(a)=\langle \pi_{i,j}(a)\xi_{i,j},\xi_{i,j}\rangle$ for any $a\in {\mathcal A}_{i,j}$.
For any $i,j$, let $\lambda_{i,j}$ be the *-representation of ${\mathcal A}_{i,j}$ given by
$$
\lambda_{i,j}(a)=V_{i,j}(\pi_{i,j}(a)\otimes I_{{\mathcal H}(i,j)}) V_{i,j}^{*} \;\;{\rm for}\;\; a\in {\mathcal A}_{i,j},
$$
where $I_{{\mathcal H}(i,j)}$ denotes the identity on ${\mathcal H}(i,j)$.
Note that these representations are, in general, non-unital.
In fact,
$$
\lambda_{i,j}(1_{i,j})=V_{i,j}V_{i,j}^{*}=r_{i,j}+s_{i,j}
$$
where $r_{i,j}$ and $s_{i,j}$ are canonical projections in $B({\mathcal H})$ given by
$$
r_{i,j}=P_{\mathcal{H}(i,j)}\;\;\;{\rm and}\;\;\; s_{i,j}=P_{\mathcal{K}(i,j)}
$$
where $\mathcal{K}(i,j)=\mathcal{H}_{i,j}^{0}\otimes \mathcal{H}(i,j)$. For given
$i,j$, the projections $r_{i,j}$ and $s_{i,j}$ are orthogonal and 
their sum is the canonical projection onto the subspace of ${\mathcal H}$ onto which $\lambda({\mathcal A}_{i,j})$ 
acts non-trivially. 

The $\lambda_{i,j}$ remind the representations $\lambda_{i}$ of free probability, 
but the corresponding operators $\lambda_{i,j}(a)$ have larger kernels. 
Using $\lambda_{i,j}$'s, we shall define product representations 
on 
$$
{\mathcal A}:=\sqcup_{i,j}{\mathcal A}_{i,j},
$$ 
the free product without identification of units, equipped with the unit $1_{\mathcal A}$,
and products of states which are analogs of the free product representation and the free product of states, respectively.
\begin{Definition}
{\rm The {\it matricially free product representation}
$\pi_{M}=*_{i,j}^{M}\pi_{i,j}$ is the unital *-homomorphism 
$\lambda: {\mathcal A}\rightarrow B({\mathcal H})$ given by the linear extension of 
$$
\lambda(1_{\mathcal A})={\bf 1}\;\;{\rm and}\;\; \lambda(a_1a_2\ldots a_n)
=\lambda_{i_1,j_1}(a_1)\lambda_{i_2,j_2}(a_2)\ldots \lambda_{i_n,j_n}(a_n)
$$ 
for any $a_{k}\in {\mathcal A}_{i_k,j_k}$, $k=1, \ldots, n$, with 
$(i_1,j_1)\neq (i_2,j_2)\neq \ldots \neq (i_n,j_n)$. 
The associated state $\varphi=*^{M}_{i,j}\varphi_{i,j}:{\mathcal A}\rightarrow {\mathbb C}$ is given by
$$
\varphi(a)=\langle \pi_{M}(a)\xi,\xi\rangle
$$ 
and will be called the {\it matricially free product of} 
$(\varphi_{i,j})$.}
\end{Definition}

Basic properties of the product state $\varphi$ are collected in the propositions given below.
Roughly speaking, they show that this state (on the free product of $C^{*}$-algebras without identification of units) has similar properties as the free product of states (on the free product of $C^{*}$-algebras with identification of units) except that the units of these algebras act as units only on `matricial' tensor products and otherwise they act as null projections.

For that purpose, it will be useful to introduce the following sets of indices:
$$
\Lambda_{n}=\{((i_1,i_2),(i_2,i_3),\ldots , (i_n,i_{n})): (i_1,i_2)\neq (i_2,i_3)\neq \ldots \neq (i_n,i_n)\}
$$
and $\Lambda=\bigcup_{n=1}^{\infty} \Lambda_{n}$. Clearly, the set 
$\Lambda$ encodes the matricial, freeness and diagonal subordination properties.
Finally, ${\mathcal I}$ stands for the unital subalgebra of ${\mathcal A}$ generated by  
the units $1_{i,j}$.

\begin{Proposition}
Let $\varphi$ be the matricially free product of states $(\varphi_{i,j})$
and let $a_{k}\in {\mathcal A}_{i_k,j_k}$, where $k \in [n]$ and
$(i_1,j_1)\neq \ldots \neq (i_n,j_n)$.  
\begin{enumerate}
\item
If $a_{k}\in {\rm Ker}\,\varphi_{i_k,j_k}$ for $k\in [n]$, then $\varphi(a_1a_2\ldots a_n)=0$.
\item

If $a_r=1_{i_r,j_r}$ and $a_{m}\in {\rm Ker}\,\varphi_{i_m,j_m}$ for $r < m \leq n$, then
$$
\varphi(a_1\ldots a_n)=\left\{
\begin{array}{cc}
\varphi(a_1\ldots a_{r-1}a_{r+1}\ldots a_n)&\;if\;\; ((i_r,j_r), \ldots , (i_n,j_n))\in \Lambda\\
0 & otherwise
\end{array}\right..
$$
\item
For any $a\in {\mathcal A}$, $u_1,u_2\in {\mathcal I}$ and $i,j\in I$, it holds that 
$$
\varphi(u_1au_2)=\varphi(u_1)\varphi(a)\varphi(u_2)\;\; and\;\; 
\varphi(1_{i,j})=\delta_{i,j}.
$$
\item
The restriction of $\varphi$ to ${\mathcal A}_{j,j}$ is $\varphi_{j,j}$ for any $j\in I$.
\item
The mixed moments $\varphi(a_1a_2\ldots a_n)$ 
are uniquely expressed in terms of mixed moments of products of $a_{k}$'s 
in the states $\varphi_{i_k,j_k}$.
\end{enumerate}
\end{Proposition}
{\it Proof.}
If $a_{k}\in {\rm Ker}\,\varphi_{i_k,j_k}$ for $k \in[n]$, where
$(i_{1},j_{1})\neq \ldots \neq (i_n,j_n)$, then it follows from the
definition of the $\lambda_{i,j}$ that
\begin{enumerate}
\item
$\pi_{M}(a_{1}\ldots a_{n})\xi=0$ if $((i_1,j_1), \ldots , (i_n,j_n))\notin \Lambda$
\item
$\pi_{M}(a_{1}\ldots a_{n})\xi\in {\mathcal H}_{i_1,j_1}^{0}\otimes \ldots \otimes {\mathcal H}_{i_n,j_n}^{0}\;$ if $\;((i_1,j_1), \ldots , (i_n,j_n))\in \Lambda$.
\end{enumerate}
In both cases we obtain a vector orthogonal to $\xi$ on the RHS, which proves (1).
Suppose now that the assumptions of (2) hold.
If $((i_r,j_r), \ldots , (i_n,j_n))\in \Lambda$, then $\lambda_{i_{r},j_{r}}(1_{i_r,j_r})$
acts as a unit on ${\mathcal H}_{i_{r+1},j_{r+1}}^{0}\otimes \ldots \otimes {\mathcal H}_{i_{n},j_{n}}^{0}$ by the definition of the representations $\lambda_{i,j}$.
On the other hand, $\lambda_{i_{r},j_{r}}(1_{i_r,j_r})$ kills
any simple tensor beginning with $h\in {\mathcal H}_{i_{r+1},j_{r+1}}^{0}$ if 
$j_{r}\neq i_{r+1}$ or $((i_{r+1},j_{r+1}), \ldots , (i_n,j_n))\notin \Lambda$
since $V_{i_r,j_r}$ does, which completes the proof of (2).
In turn, (3) follows from the action of the $\lambda(1_{i,j})$ onto $\xi$.
That $\varphi$ agrees with $\varphi_{j,j}$ on
${\mathcal A}_{j,j}$ for any $j\in I$ follows from the action
of the $\lambda_{j,j}(a)$, $a\in {\mathcal A}_{j,j}$, onto $\xi$, namely
$\pi_{M}(a)\xi=(\pi_{j,j}(a)\xi)^{0}+ \varphi_{j,j}(a)\xi$,
which gives (4). 
Finally, (5) is a consequence of (1)-(2).
\hfill $\blacksquare$\\

In a similar way we can define states associated with other unit vectors from ${\mathcal H}$.
For our purposes, we will need states associated with unit vectors 
$e_{j}\in {\mathcal H}_{j,j}^{0}$ which are in the ranges of $\pi_{j,j}({\mathcal A}_{j,j})$, where 
$j\in I$, respectively, namely $\varphi_{j}:{\mathcal A}\rightarrow {\mathbb C}$
defined by the formulas
$$
\varphi_{j}(a)=\langle {\pi}_{M}(a)e_{j},e_{j}\rangle.
$$
These states will be called {\it conditions} associated with $\varphi$
and will be used for computing traces.
They have similar properties as $\varphi$ as the proposition given 
below shows in more detail. The main difference is that $\varphi_j$ 
satisfies the condition of freeness type only for indices which satisfy 
$(j,j)\neq (i_1,j_1)\neq \ldots \neq (i_n,j_n)\neq (j,j)$.
Moreover, $\varphi_{j}|{\mathcal I}$ is quite different than $\varphi|{\mathcal I}$ due to
different normalization conditions. 
\begin{Proposition}
Let $\varphi_j$, where $j\in I$, be the conditions associated with $\varphi=*^{M}_{i,j}\varphi_{i,j}$
and let  let $a_{k}\in {\mathcal A}_{i_k,j_k}$, where $k \in [n]$ and
$(j,j)\neq (i_1,j_1)\neq \ldots \neq (i_n,j_n)\neq (j,j)$. Then
\begin{enumerate}
\item
If $a_{k}\in {\rm Ker}\,\varphi_{i_k,j_k}$ for $k \in [n]$, then
$\varphi_{j}(a_1a_2\ldots a_n)=0$ for each $j$.
\item
If $a_r=1_{i_r,j_r}$ and $a_{m}\in {\rm Ker}\,\varphi_{i_m,j_m}$ for $r < m \leq n$, then
$$
\varphi_{j}(a_1\ldots a_n)=\left\{
\begin{array}{cc}
\varphi_{j}(a_1\ldots a_{r-1}a_{r+1}\ldots a_n)&\;if\;\; ((i_r,j_r), \ldots , (i_n,j_n))\in \Lambda\\
0 & otherwise
\end{array}\right..
$$
\item
For any $a\in {\mathcal A}$, $u_1,u_2\in {\mathcal I}$ and $i,j,k\in I$, it holds that 
$$
\varphi_{j}(u_1au_2)=
\varphi_{j}(u_1)\varphi_{j}(a)\varphi_{j}(u_2)\;\;and\;\;
\varphi_{j}(1_{i,k})=\delta_{j,k}.
$$
\item
The restriction of $\varphi_{j}$ to ${\mathcal A}_{i,j}$ is $\varphi_{i,j}$ for $i\neq j$.
\item
The mixed moments $\varphi_{j}(a_1a_2\ldots a_n)$ 
are uniquely expressed in terms of mixed moments of products of $a_{k}$'s 
in the states $\varphi_{i_k,j_k}$.
\end{enumerate}
\end{Proposition}
{\it Proof.}
The proof is similar to that of Proposition 2.2.
The main difference concerns (3) and (4). In this context, notice that
the unit vectors $e_j$ play the same role with respect to the action of 
the $\lambda_{i,j}(a)$ for any $i\neq j$ as $\xi_{i,j}$ plays with respect to the action of $\pi_{i,j}(a)$, 
where $a\in {\mathcal A}_{i,j}$, and thus
$\varphi_{j}(a)=\langle \lambda_{i,j}(a)e_j,e_j\rangle =\langle \pi_{i,j}(a)\xi_{i,j},\xi_{i,j}\rangle 
=\varphi_{i,j}(a)$, which gives (4).
\hfill $\blacksquare$
\begin{Remark}
{\rm Note that states $\varphi$ and $(\varphi_{j})$ share together the property of extending the array of states $(\varphi_{i,j})$.
Thus, $\varphi$ extends the diagonal states $\varphi_{j,j}$ for all $j$, 
but it does not extend the non-diagonal states $\varphi_{i,j}$ for $i\neq j$ 
since $\pi_{M}(a)\xi=0$ for any $a\in {\mathcal A}_{i,j}$.
In turn, $\varphi_{j}$ extends $\varphi_{i,j}$ for $i\neq j$, but it does not
extend $\varphi_{j,j}$ or any of the other diagonal states. 
This is a natural consequence of differences in
the definitions of the diagonal and non-diagonal partial isometries.
Moreover, it is clear that for each 
$\varphi_j$ there exists $b_{j}\in {\mathcal A}_{j,j}\cap {\rm Ker}\varphi$ such that
$$
\varphi_{j}(w)=\varphi(b_j^{*}wb_j)
$$
for any $w\in \sqcup_{i,j}\mathcal{A}_{i,j}$.
Therefore we can reduce computations of mixed moments in the state $\varphi_j$ 
to computations of mixed moments in the state $\varphi$.}
\end{Remark} 

Finally, let us denote by $\lambda({\mathcal I})$ 
the unital commutative *-subalgebra 
of $B({\mathcal H})$ generated by the 
$\lambda(1_{i,j})$, where $i,j\in I$. 
By abuse of notation, $\lambda(1_{i,j})$ will also
be denoted by $1_{i,j}$. Observe that these projections are not 
mutually orthogonal and that is why it is often convenient to use their subprojections $r_{i,j}$ and $s_{i,j}$.

\section{Strongly matricially free products} 
Of special importance is the subspace of 
the matricially free Fock space, called the `strongly matricially free Fock space',
in which the diagonal Hilbert spaces appear only at the end of tensor products. 
The main reason is that it is related to both free and monotone Fock spaces.
We also study the associated product states which can be viewed as direct generalizations of both 
free and monotone products of states. 
\begin{Definition}
{\rm By the {\it strongly matricially free Fock space} over 
$\widehat{{\mathcal H}}:=({\mathcal H}_{i,j})$
we understand the subspace of $\mathcal{M}(\widehat{\mathcal H})$ of the form
$$
\mathcal{R}(\widehat{\mathcal H})=
{\mathbb C}\Omega \oplus \bigoplus _{m=1}^{\infty}
\bigoplus_{\stackrel {i_{1}\neq \ldots \neq i_{m}}
{\scriptscriptstyle n_1, \ldots , n_m \in {\mathbb N}}}
{\mathcal H}_{i_1,i_2}^{\otimes \,n_1}
\otimes 
{\mathcal H}_{i_2,i_3}^{\otimes \,n_2}
\otimes \ldots \otimes 
{\mathcal H}_{i_{m},i_{m}}^{\otimes \,n_m}.
$$
with the canonical inner product.}
\end{Definition}

A justification for the word `strong' is that in this case the words $i_1i_2\ldots i_m$ which label 
the tensor products in the above definition satisfy $i_1\neq i_2\neq \ldots \neq i_m$.
\begin{Example}
{\rm The simplest space of this type is associated with a 2-dimensional square array 
$\widehat{{\mathcal H}}$. Then 
$$
\mathcal{R}(\widehat{\mathcal H})
=\bigoplus_{m=0}^{\infty}\mathcal{R}^{(m)}(\widehat{\mathcal H}),
$$
where the first few summands are of the form
\begin{eqnarray*}
\mathcal{R}^{(0)}(\widehat{\mathcal H})&=&{\mathbb C}\Omega\\
\mathcal{R}^{(1)}(\widehat{\mathcal H})&=&{\mathcal H}_{1,1}
\oplus {\mathcal H}_{2,2}\\
\mathcal{R}^{(2)}(\widehat{\mathcal H}) &=&
{\mathcal H}_{1,1}^{\otimes 2}\oplus 
{\mathcal H}_{2,2}^{\otimes 2}\oplus
({\mathcal H}_{1,2}\otimes {\mathcal H}_{2,2})\oplus
({\mathcal H}_{2,1}\otimes {\mathcal H}_{1,1}) \\
\mathcal{R}^{(3)}(\widehat{\mathcal H})&=&
{\mathcal H}_{1,1}^{\otimes 3}
\oplus {\mathcal H}_{2,2}^{\otimes 3}
\oplus 
({\mathcal H}_{2,1}\otimes {\mathcal H}_{1,1}^{\otimes 2})
\oplus 
({\mathcal H}_{1,2}\otimes {\mathcal H}_{2,2}^{\otimes 2})
\oplus 
({\mathcal H}_{2,1}^{\otimes 2}\otimes {\mathcal H}_{1,1})\\
&&
\oplus 
({\mathcal H}_{1,2}^{\otimes 2}\otimes {\mathcal H}_{2,2})
\oplus 
({\mathcal H}_{1,2}\otimes {\mathcal H}_{2,1}\otimes {\mathcal H}_{1,1})
\oplus 
({\mathcal H}_{2,1}\otimes {\mathcal H}_{1,2}\otimes {\mathcal H}_{2,2}),
\end{eqnarray*}
etc. In contrast to
$\mathcal{M}(\widehat{\mathcal H})$, we do not have tensor products like
${\mathcal H}_{2,2}\otimes {\mathcal H}_{2,1}\otimes {\mathcal H}_{1,1}$ 
and ${\mathcal H}_{1,1}\otimes {\mathcal H}_{1,2}\otimes {\mathcal H}_{2,2}$
in the summand of the third order.
}
\end{Example}
\begin{Remark}
{\rm For a given array of Hilbert spaces $\widehat{\mathcal H}=({\mathcal H}_{i,j})$, we have
inclusions
$$
{\mathcal R}(\widehat{\mathcal H})\subseteq {\mathcal M}(\widehat{\mathcal H})\subseteq 
{\mathcal F}(\bigoplus_{i,j}{\mathcal H}_{i,j})
$$
which, in most cases, are proper.
Moreover, if we have a square array
and ${\mathcal H}_{i,j}\cong {\mathcal H}_{i}$ for any $i,j\in I$, 
where $({\mathcal H}_{i})_{i\in I}$ is a family of Hilbert spaces, then there is 
a natural isomorphism
$$
\mathcal{R}(\widehat{\mathcal H})\cong \mathcal{F}(\bigoplus_{i\in I}{\mathcal H}_{i})
$$
since ${\mathcal H}_{i_1,i_2}^{\otimes n_1}\otimes {\mathcal H}_{i_2,i_3}^{\otimes n_2}\otimes \ldots \otimes {\mathcal H}_{i_m,i_m}^{\otimes n_m}\cong {\mathcal H}_{i_1}^{\otimes n_1}\otimes {\mathcal H}_{i_2}^{\otimes n_2}\otimes \ldots \otimes {\mathcal H}_{i_m}^{\otimes n_m}$ for any $i_1,i_2,\ldots, i_m\in I$, $n_1,\ldots n_m, m\in {\mathbb N}$.
Similarly, if we have a lower-triangular array and ${\mathcal H}_{i,j}\cong {\mathcal H}_{i}$ for any $i\geq j$, then 
$\mathcal{R}(\widehat{\mathcal H})$ is isomorphic to the monotone Fock space. }
\end{Remark}

Moreover, as expected, there is a product of Hilbert spaces related
to the strongly matricially free Fock space, and an analog of Proposition 2.1 holds.

\begin{Definition}
{\rm By the {\it strongly matricially free product 
of $({\mathcal H}_{i,j},\xi_{i,j})$} we understand 
the pair $({\mathcal G}, \xi)$, where ${\mathcal G}$ 
is the subspace of ${\mathcal H}$ of the form
$$
{\mathcal G}= {\mathbb C}\xi \oplus\bigoplus_{m=1}^{\infty}
\bigoplus_{i_{1}\neq \ldots \neq i_{m}}
{\mathcal H}_{i_1,i_2}^{0}
\otimes 
{\mathcal H}_{i_2,i_3}^{0}
\otimes \ldots \otimes 
{\mathcal H}_{i_{m},i_{m}}^{0}
$$
We denote it $({\mathcal G}, \xi)=*^{R}_{i,j}({\mathcal H}_{i,j},\xi_{i,j})$.}
\end{Definition}

If we consider a family of unital $C^{*}$-algebras $({\mathcal A}_{i})_{i\in I}$, each equipped with 
a family of states $(\varphi_{i,j})_{j\in I}$, then we can look at this product space
as follows. If $({\mathcal H}_{i,j},\pi_{i,j},\xi_{i,j})$ is the GNS triple associated with 
the pair $({\mathcal A}_{i}, \varphi_{i,j})$, then the 
Hilbert space $\mathcal{H}_{i,j}$ (as well as the corresponding state and representation) 
taken as the representation space for the algebra ${\mathcal A}_{i}$ at some given tensor site
depends on the algebra $\mathcal{A}_{j}$ represented at the 
following tensor site on the space ${\mathcal H}_{j,k}$ for some $k$.
It is worth noting that in this framework we can assume that for fixed $i\in I$ all vectors $\xi_{i,j}, j\in I$, are identified since we can take the tensor product of Hilbert spaces $\otimes_{j}{\mathcal H}_{i,j}$ and set $\xi_{j}=\otimes_{j\in I}\xi_{i,j}$ for each $i\in I$. It is not hard to see that in this framework our model 
is related to freeness with infinitely many states [8,9]. 

The construction of the product state is similar to that in the usual case.
The only difference in all definitions 
is that the sets $\Lambda_{n}$ are replaced by 
$$
\Gamma_{n}=\{((i_1,i_2),(i_2,i_3),\ldots , (i_n,i_{n})): i_1\neq i_2\neq \ldots \neq i_n\}
$$ 
and their union $\Lambda$ by $\Gamma=\bigcup_{n=1}^{\infty} \Gamma_{n}$.
Note that the conditions which define the $\Gamma_{n}$ are stronger than those
which define the $\Lambda_n$ and therefore all objects constructed in the strong
case are obtained from the standard ones by a projection-type operation. 
 
The partial isometries in the strong case, denoted by $W_{i,j}$, remind the $V_{i,j}$ except that they refer 
to ${\mathcal G}$ rather than ${\mathcal H}$. In particular, the diagonal partial isometries
$W_{j,j}:{\mathcal H}_{j,j}\otimes {\mathcal G}(j,j)\rightarrow {\mathcal G}$ are given by
$$
\xi_{j,j}\otimes \xi \rightarrow \xi\;\;\;{\rm and}\;\;\;
{\mathcal H}_{j,j}^{0}\otimes \xi \rightarrow  {\mathcal H}_{j,j}^{0}
$$
where ${\mathcal G}(j,j)={\mathbb C}\xi$ for any $j$, whereas the non-diagonal partial 
isometries $W_{i,j}$ and the associated subspaces ${\mathcal G}(i,j)$ are similar 
to the ${\mathcal H}(i,j)$.

\begin{Definition}
{\rm Let $\rho_{i,j}$ be the *-representation of ${\mathcal A}_{i,j}$ defined as in the
standard case, with $V_{i,j}$ replaced by $W_{i,j}$ for any $i,j$. The corresponding 
{\it strongly matricially free product representation} $\pi_{R}=*^{R}_{i,j}\pi_{i,j}$ and 
{\it strongly matricially free product of states} $*_{i,j}^{R}\varphi_{i,j}$ are 
defined in terms of the $\rho_{i,j}$ as in the standard case.}
\end{Definition}
The basic results on the strongly matricially free product of states 
are similar to those in Propositions 2.2-2.4 and therefore we state 
them in an abbreviated form without a proof.
\begin{Proposition}
Let $\varphi$ be the strongly matricially free product of states $(\varphi_{i,j})$ and let
$\varphi_{j}$ be the state associated with any unit vector from ${\mathcal H}_{j,j}^{0}$
which is in the range of $\rho_{j,j}({\mathcal A}_{j,j})$, $j\in I$.
Then the statements of Propositions 2.2-2.3 remain true, with $\Lambda$ replaced by $\Gamma$.
\end{Proposition}
As we mentioned earlier, one of the advantages of using the strong structures
is that they are straightforward generalizations of those in free probability
(in the case of square arrays) and monotone probability (in the case of triangular arrays).
\begin{Remark}
{\rm 
Using the strongly matricially free product of states, we can obtain the free product of 
states if we take a square array. Namely, if ${\mathcal A}_{i,j}={\mathcal A}_{i}$
and $\varphi_{i,j}=\varphi_{i}$ for all $i,j\in I$, then, on the level of Hilbert spaces, we 
have ${\mathcal H}_{i,j}={\mathcal H}_{i}$ and $\xi_{i,j}=\xi_{i}$ for all $i,j\in I$. 
Thus 
$$
({\mathcal G}, \xi)\cong *_{i\in I}({\mathcal H}_{i},\xi_{i})
$$ 
and $\rho_{i}(a):=\sum_{j}\rho_{i,j}(a)$, $a\in {\mathcal A}_{i}$, understood as the strong limit, 
extends to a (unital) *-representation of 
${\mathcal A}_{i}$ on ${\mathcal H}$. Then $\rho=\sqcup_{i\in I}\rho_{i}$ agrees with the 
free product representation of $\sqcup_{i\in I}{\mathcal A}_{i}$ on ${\mathcal G}$
and therefore the corresponding strongly matricially free product of states gives the free product of states.
If $\varphi_{i,i}=\varphi_i$ and $\varphi_{i,j}=\psi_{i}$ for any $i\neq j$, we obtain
the conditionally free product of states.}
\end{Remark}
\begin{Remark}
{\rm 
The monotone product of states [17] is obtained from the strongly matricially free product of states
when we take a triangular array $T$, with
${\mathcal A}_{i,j}={\mathcal A}_{i}$ and $\varphi_{i,j}=\varphi_{i}$ 
for all $(i,j)\in T\subset I\times I$ 
(it can be obtained from the square array by setting ${\mathcal A}_{i,j}=0$ for $i<j$), where
$I$ is a totally ordered set. 
Then ${\mathcal H}_{i,j}={\mathcal H}_{i}$ and $\xi_{i,j}=\xi_{i}$ for $i\geq j$ 
and 
$$
({\mathcal G}, \xi)\cong \vartriangleright_{i\in I}({\mathcal H}_{i},\xi_{i}),
$$ 
the expression on the right-hand side being the monotone product of Hilbert spaces. 
Moreover, $\tau_{i}(a):=\sum_{j \leq i}\rho_{i,j}(a)$ extends 
to a (non-unital) *-representation of ${\mathcal A}_{i}$ on ${\mathcal G}$, such that
$\tau=\sqcup_{i\in I}\tau_{i}$ agrees with the 
monotone product representation of $\sqcup_{i\in I}{\mathcal A}_{i}$ on ${\mathcal G}$
and therefore the corresponding strongly matricially free product of states gives the monotone product of states.} 
\end{Remark}

\section{Matricial freeness}
Guided by the notion of the matricially free product of states, we shall introduce now the
associated concept of independence called `matricial freeness', and closely related to it, 
`strong matricial freeness'. They involve 
arrays of noncommutative probability spaces and to some extent they remind models with many states [6,8,9], 
but they cannot be reduced in a natural way to any of these (freeness with infinitely many states has some non-empty intersection with `strong matricial freeness' and conditional freeness is its special case). Moreover, we will study discrete matricially free Fock spaces, both standard and strong.

Let ${\mathcal A}$ be a unital algebra with an array $({\mathcal A}_{i,j})$ of not necessarily unital subalgebras of 
${\mathcal A}$. Assume that each ${\mathcal A}_{i,j}$ has an internal unit $1_{i,j}$ and
that the subalgebra ${\mathcal I}$ of ${\mathcal A}$ generated by all internal units is commutative.  
Further, let $\varphi$ be a distinguished state on ${\mathcal A}$ such that $\varphi(1_{i,j})=\delta_{i,j}$
for any $i,j$ and let $\{\varphi_{j}:i\in I\}$ be a family of additional states on 
${\mathcal A}$ such that $\varphi_{j}(1_{i,k})=\delta_{j,k}$ for any $i,j,k$, which will be
called {\it conditions}. Here, by a state on
${\mathcal A}$ we understand a normalized linear functional (if ${\mathcal A}$ is a *-algebra,  
we require this functional to be positive). 

In this situation it is convenient to form an array $(\varphi_{i,j})$ of states
on ${\mathcal A}$ by 
$$
\varphi_{j,j}=\varphi \;\;\;{\rm and}\;\;\;\varphi_{i,j}=\varphi_{j}\;\;{\rm for}\;\;i\neq j
$$ 
which will be said to be {\it defined} by the state $\varphi$ and the conditions $\varphi_{j}$.
\begin{Definition}
{\rm Let $(\varphi_{i,j})$ be defined by the state $\varphi$ and the conditions $\varphi_j$.
We say that $(1_{i,j})$ is a {\it matricially free array of units} 
associated with $({\mathcal A}_{i,j})$ and  $(\varphi_{i,j})$ if
\begin{enumerate}
\item
$\varphi (u_1au_2)=\varphi(u_1)\varphi(a)\varphi(u_2)$ 
for any $a\in {\mathcal A}$ and $u_1,u_2\in {\mathcal I}$,
\item
if $a_{k}\in {\mathcal A}_{i_k,j_k}\cap {\rm Ker}\varphi_{i_k,j_k}$, where $r<k\leq n$ and $r<n$, then
$$
\varphi(a1_{i_r,j_r}a_{r+1}\ldots a_n)=
\left\{
\begin{array}{cc}
\varphi(aa_{r+1} \ldots a_n) & {\rm if}\;((i_{r},j_{r}), \ldots , (i_{n},j_n))\in \Lambda\\
0 & {\rm otherwise}
\end{array}
\right..
$$
where $a\in {\mathcal A}$ is arbitrary and $(i_r,j_r)\neq \ldots \neq (i_n,j_n)$.
\end{enumerate}
The array $(1_{i,j})$ is called a {\it strongly matricially free array of units}
if $\Lambda$ is replaced by $\Gamma$.}
\end{Definition}
The above definition enables us to define the concepts of {\it matricial freeness}
and its strong version called {\it strong matricial freeness}. 
They both bear some resemblance to freeness, but the main difference is that the 
identified unit in the context of freeness is replaced by the (strongly) matricially 
free array of units. In fact, we will see later that it is the strong matricial freeness 
which can be viewed as a direct generalization of freeness.
\begin{Definition}
{\rm We say that $({\mathcal A}_{i,j})$ is 
{\it matricially free} with respect to $(\varphi_{i,j})$ if  
\begin{enumerate}
\item for any $a_{k}\in {\rm Ker}\varphi_{i_k,j_k}\cap {\mathcal A}_{i_k,j_k}$, where $k\in [n]$
and $(i_1,j_1)\neq \ldots \neq (i_n,j_n)$,
$$
\varphi(a_1a_2\ldots a_n)=0
$$
\item
$(1_{i,j})$ is a matricially free array of units associated with 
$({\mathcal A}_{i,j})$ and $(\varphi_{i,j})$.
\end{enumerate}
In an analogous manner we define {\it strongly matricially free} arrays of subalgebras.}
\end{Definition}
\begin{Definition}
{\rm 
The array of variables $(a_{i,j})$ in a unital algebra ${\mathcal A}$ will be called 
({\it strongly}) {\it matricially free} with respect to the array $(\varphi_{i,j})$ 
defined by a state $\varphi$ and the conditions $\varphi_j$ if there 
exists an array of (strongly) matricially free array of units 
$(1_{i,j})$ in ${\mathcal A}$ such that the array of algebras $({\mathbb C}[a_{i,j},1_{i,j}])$ 
is (strongly) matricially free with respect to $(\varphi_{i,j})$.}
\end{Definition}

Using the above definitions, we can compute mixed moments 
of arbitrary matricially free random variables. Namely,
writing $a_{k}=a_{k}^{0}+\varphi_{i_k,j_k}(a_{k})1_{i_k,j_k}$ 
for each $k\in [n]$, we obtain a recursion
\begin{eqnarray*}
\varphi(a_1\ldots a_n)&=&
\;\,\sum_{1\leq k \leq n}\varphi_{i_k,j_k}(a_{k})\varphi(a_1^{0}\ldots 1_{i_k,j_k}\ldots a_{n}^{0})\\
&+&
\sum_{1\leq k<l\leq n}\varphi_{i_k,j_k}(a_k)\varphi_{i_l,j_l}(a_l)
\varphi(a_1^{0}\ldots 1_{i_k,j_k}\ldots 1_{i_l,j_l}\ldots a_{n}^{0}) \\
&+& \;\;\;\;\ldots \\
&+& \;\;\;\;\varphi_{i_1,j_1}(a_1)\ldots \varphi_{i_n,j_n}(a_n)\varphi(1_{i_1,j_1}\ldots 1_{i_n,j_n})
\end{eqnarray*}
for any $(i_1,j_1)\neq \ldots \neq (i_n,j_n)$. 
It is easy to see that that the mixed moments on the right-hand side,
written here in a slightly simplified manner, reduce to moments of orders smaller than $n$.
Note that these depend on $\varphi|{\mathcal I}$ in an essential way. 
That is why the states $\varphi$ and the $\varphi_{j}$ of Section 2 give 
different mixed moments, although they have similar properties.    
\begin{Remark}
{\rm    
If ${\mathcal A}$ is a unital *-algebra, then, in addition, we require that 
the functionals $\varphi_{i,j}$ are positive, the ${\mathcal A}_{i,j}$ 
are *-subalgebras and the $1_{i,j}$ are projections. Then an array of variables $(a_{i,j})$
will be called *-({\it strongly}) {\it matricially free} if the array of *-algebras 
$({\mathbb C}\langle a_{i,j},a_{i,j}^{*},1_{i,j}\rangle )$ 
is (strongly) matricially free.}
\end{Remark}

\begin{Example}
{\rm 
Let us assume that we have a two-dimensional square array of 
variables, namely $a_{1,1}=a, a_{1,2}=a', a_{2,1}=b', a_{2,2}=b$,
which are strongly matricially free with respect to the array $(\varphi_{i,j})$ 
defined by $\varphi$ and the pair $(\varphi_1,\varphi_2)$. Using the above recursion
and then Definitions 4.1-4.2, we obtain
\begin{eqnarray*}
\varphi(ab'ab)&=&\varphi(b)\varphi_1(b')(\varphi(a^2)-\varphi^2(a)),\\
\varphi(abab)&=&\varphi^2(a)\varphi^2(b),\\
\varphi(aba'b)&=&\varphi(a)\varphi_2(a')(\varphi(b^2)-\varphi^2(b)).
\end{eqnarray*}
Note that if we assume that $\varphi_1=\varphi_2=\psi$, then the sum of these
moments gives $\varphi(abab)$ for $a,b$ conditionally independent with respect to
$(\varphi, \psi)$.
In fact, the above sum is equal to $\varphi((a+a')(b+b')(a+a')(b+b'))$.
In turn, if we take the lower-triangular subarray
and take $\varphi_1$ and $\varphi$--distributions to be equal,   
the sum of the first two moments (the third one is zero since there is no $a'$) 
gives $\varphi(abab)$ for $a,b$ monotone independent with respect to $\varphi$ 
Similar agreements hold for mixed moments of higher orders, which is in agreement
with Remarks 3.2-3.3. }
\end{Example}

\begin{Definition}
{\rm By a {\it discrete matricially free Fock space} we understand 
${\mathcal M}={\mathcal M}(\widehat{\mathcal H})$, where
${\mathcal H}_{i,j}={\mathbb C}e_{i,j}$ for any 
$(i,j)\in {\mathbb N}$, and the array $(e_{i,j})$ forms an orthonormal basis of some Hilbert space.
In an analogous manner we define the {\it discrete strongly matricially free Fock space} 
${\mathcal R}={\mathcal R}(\widehat{\mathcal H})$.}
\end{Definition}
Both ${\mathcal M}$ and ${\mathcal R}$ are subspaces of the discrete free 
Fock space ${\mathcal F}(\bigoplus_{i,j}{\mathbb C}e_{i,j})$
and thus allow for the canonical action of free creation and annihilation operators.
In order to specify the arrays of units, let us distinguish two types of their subspaces.

In the case of ${\mathcal M}$ these are 
\begin{enumerate}
\item
${\mathcal M}(i,j)$, spanned be simple tensors which begin with $e_{j,k}$ for some $k$, 
where $(j,k)\neq (i,j)$, and, in addition, by $\Omega$ if $i=j$,  
\item
${\mathcal K}(i,j)$, spanned by vectors which begin with $e_{i,j}$, where $i,j$ are arbitrary,
\end{enumerate}
and the direct sum ${\mathcal K}(i,j)\oplus {\mathcal M}(i,j)$ is the subspace of $\mathcal{M}$ 
onto which the *-algebra generated by $\ell(e_{i,j})$ restricted to ${\mathcal M}$ acts non-trivially, 
where $\ell(e_{i,j})$ denotes the canonical free creation operator
associated with vector $e_{i,j}$. 

The canonical projections onto such 
direct sums are natural candidates for the matricially free units $1_{i,j}$ and therefore we set 
$$
1_{i,j}:=P_{{\mathcal K}(i,j)\oplus {\mathcal M}(i,j)}\;\;\;{\rm and}\;\;\;
\ell_{i,j}=\ell(e_{i,j})1_{i,j}
$$
for any $i,j$. The adjoint of $\ell_{i,j}$ will be denoted by $\ell_{i,j}^{*}$.
Note that the projection $1_{i,j}$ is an internal unit in the *-algebra 
${\mathcal A}_{i,j}={\mathbb C}\langle \ell_{i,j}, \ell_{i,j}^{*} \rangle$ 
and $\ell_{i,j}^{*}\ell_{i,j}=1_{i,j}$ for any $i,j$. 
However, in the algebra ${\mathbb C}\langle \ell_{i,j},\ell_{i,j}^{*}: i,j\in {\mathbb N} \rangle$
there are more relations as the proposition given below demonstrates.
\begin{Proposition} 
In the algebra ${\mathbb C}\langle \ell_{i,j},\ell_{i,j}^{*}: i,j\in {\mathbb N} \rangle$
the following relations hold:
\begin{enumerate}
\item 
$\ell_{i,j}^{*}\ell_{i,j}=1_{i,j}$ for any $i,j$,
\item
$\ell_{i,j}^{*}\ell_{k,l}=0$, $\ell_{i,j}\ell_{k,l}=0$ and $1_{i,j}\ell_{k,l}=0$ whenever $(i,j)\neq (k,l)$ and $j\neq k$, 
\item
$1_{i,j}\ell_{j,k}=\ell_{j,k}$ for any $j,k$.
\end{enumerate}
\end{Proposition}
{\it Proof.}
We omit the elementary proof.
\hfill $\blacksquare$\\

In the case of ${\mathcal R}$, we proceed in a completely analogous fashion 
and distinguish the following subspaces:
$$
{\mathcal R}(i,j)={\mathcal M}(i,j)\cap {\mathcal R}\;\;\;{\rm and}\;\;\;
{\mathcal L}(i,j)={\mathcal K}(i,j)\cap {\mathcal R}
$$
which lead to the definitions of the strongly matricially free array of units and of the 
creation operators, respectively,
$$
1_{i,j}=P_{{\mathcal L}(i,j)\oplus {\mathcal R}(i,j)}\;\;\;{\rm and}\;\;\;
k_{i,j}=\ell(e_{i,j})1_{i,j}
$$
where, slightly abusing notation, we use the same symbols for the units as before.

Finally, we specify the states: $\varphi$ will denote the
vacuum state associated with $\Omega$, the $\varphi_{j}$ will be the states associated with
the $e_{j,j}$, $j\in {\mathbb N}$, and $(\varphi_{i,j})$ will be the array defined by $\varphi$ 
and the $\varphi_{j}$. Here, the same notation is used for states on $B({\mathcal M})$ and
$B({\mathcal R})$.
\begin{Proposition}
The array of *-subalgebras ${\mathcal A}_{i,j}={\mathbb C}\langle \ell_{i,j}, \ell_{i,j}^{*}\rangle$ of $B({\mathcal M})$, where $i,j\in {\mathbb N}$, is matricially free with respect to $(\varphi_{i,j})$. 
The array of *-subalgebras 
${\mathcal B}_{i,j}={\mathbb C}\langle k_{i,j}, k_{i,j}^{*}\rangle$ 
of $B({\mathcal R})$, where $i,j\in {\mathbb N}$, is strongly matricially free with 
respect to $(\varphi_{i,j})$. 
\end{Proposition}
{\it Proof.}
The proof is similar to that of Voiculescu for the discrete free Fock space
given in [21]. We shall look at the case of $(\mathcal{A}_{i,j})$ since
the case of $(\mathcal{B}_{i,j})$ is analogous. Each algebra $\mathcal{A}_{i,j}$ is spanned by 
operators of the form
$$
\ell_{i,j}^{q}\ell_{i,j}^{*p},\;\;{\rm where}\;\; p+q>0,
$$
and the projection $1_{i,j}$. However, the corresponding moments vanish:
$$
\varphi_{i,j}(\ell_{i,j}^{q}\ell_{i,j}^{*p})=0
$$
for any $i,j$ since $p+q>0$. Moreover,
$\varphi_{i,j}(1_{i,j})=1$ for any $i,j$.
Therefore, in order to show that condition (1) of Definition 4.2 holds, 
it is enough to show that
$$
\varphi(\ell_{i_1,j_1}^{q_1}\ell_{i_1,j_1}^{*p_1}\ldots 
\ell_{i_n,j_n}^{q_n}\ell_{i_n,j_n}^{*p_n})=0
$$
whenever $(i_1,j_1)\neq \ldots \neq (i_n,j_n)$ and $p_1+q_1>0, \ldots, p_n+q_n>0$.
The same argument as in [21] allows us to reduce the proof to the case when $q_1=\ldots =q_n=0$,
which implies that $p_1>0, \ldots , p_n>0$. But then the moment clearly vanishes. This
proves condition (1) of Definition 4.2. Condition (2) follows easily from the definition of the 
projections $1_{i,j}$ in view of the relations given in Proposition 4.1.  This completes the proof.
\hfill$\blacksquare$
\begin{Example}
{\rm For $i,j\in I$, let $G_{i,j}\cong F(1)$ be the free group 
on one generator $g_{i,j}$ with unit $\epsilon_{i,j}$, and
let $\lambda_{i,j}$ be the corresponding unitary operator on 
the space $l^{2}(G_{i,j})$ given by $\lambda_{i,j}(g)\delta (h)=\delta (gh)$,
where $\{\delta (g): g\in G_{i,j}\}$ is the canonical basis. Consider the subspace
$l^{2}_{M}$ of $l^{2}(*_{i,j}G_{i,j})$ spanned by vectors of the form $\delta (g)$, where
$g$ is either the unit $e$ of the free product $*_{i,j}G_{i,j}$, or 
a product of the form $g_1g_2\ldots g_{m}$, where 
$g_{k}\in G_{i_k,i_{k+1}}^{0}:=G_{i_{k} i_{k+1}}\setminus \{\epsilon_{i_{k},i_{k+1}}\}$ for each $k$, with $(i_1,i_2)\neq \ldots \neq (i_{m},i_{m+1})$ and $i_{m+1}=i_m$. The space $l_{M}^{2}$ can be viewed as the space of
square integrable functions on the (non-existent) `matricially free product of groups'.
Let $1_{i,j}$ denote the projection from $l^{2}(*_{i,j}G_{i,j})$ onto the subspace
of $l^{2}_{M}$ spanned by vectors $\delta (g)$, where $g$ begins with an element
from $G_{i,j}^{0}$ or $G_{j,k}^{0}$ for some $k$ or, in the case of $i=j$, also $g=e$.
Define
$$
\widetilde{\lambda}_{i,j}(g)=\lambda_{i,j}(g)1_{i,j}
$$
for $g\in G_{i,j}$. Then the array $({\mathcal A}_{i,j})$, where
${\mathcal A}_{i,j}$ is the *-subalgebra of $B(l^2(*_{i,j}G_{i,j}))$ generated by $\widetilde{\lambda}_{i,j}(g_{i,j})$ and $1_{i,j}$,
with the standard involution, is matricially free with respect to the array
$(\varphi_{i,j})$, where the diagonal states coincide with
$\varphi(.)=\langle .\delta (e), \delta (e)\rangle$ and the non-diagonal ones
in the $j$-th column coincide with 
$\varphi_{j}(.)=\langle .\delta (g_{j,j}), \delta (g_{j,j})\rangle$.}
\end{Example}
\begin{figure}
\unitlength=1mm
\special{em.linewidth 0.5pt}
\linethickness{0.5pt}
\begin{picture}(100.00,70.00)(-25.00,-10.00)
%%%%%%%horizontal lines%%%%%%%%%%%%%%%%%%%%%%%%%%%%%%%
%%-2 horizontal line 
\put(15.00,0.00){\line(1,0){20.00}}
%%its vertices
\put(25.00,0.00){\circle*{1.50}}
\put(20.00,0.00){\circle{1.50}}
\put(22.00,-2.00){\circle*{1.50}}
\put(18.00,2.00){\circle*{1.50}}
\put(22.00,2.00){\circle{1.50}}
\put(18.00,-2.00){\circle{1.50}}
\put(30.00,0.00){\circle{1.50}}
\put(28.00,-2.00){\circle*{1.50}}
\put(32.00,2.00){\circle*{1.50}}
\put(28.00,2.00){\circle{1.50}}
\put(32.00,-2.00){\circle{1.50}}
\put(16.50,0.00){\circle{1.50}}
\put(33.50,0.00){\circle{1.50}}
%%with stars
\put(16.50,-3.50){\line(1,1){7.00}}
\put(23.50,-3.50){\line(-1,1){7.00}}
\put(26.50,-3.50){\line(1,1){7.00}}
\put(33.50,-3.50){\line(-1,1){7.00}}
%%%%%%%%%%%%%%%%%%%%%%%%%%%%%
%-1 horizontal line
\put(13.00,10.00){\line(1,0){24.00}}
%%its vertices
\put(32.50,10.00){\circle{1.50}}
\put(17.50,10.00){\circle{1.50}}
\put(15.50,12.00){\circle*{1.50}}
\put(15.50,8.00){\circle{1.50}}
\put(19.50,12.00){\circle{1.50}}
\put(19.50,8.00){\circle*{1.50}}
\put(30.50,12.00){\circle{1.50}}
\put(30.50,8.00){\circle*{1.50}}
\put(34.50,12.00){\circle*{1.50}}
\put(34.50,8.00){\circle{1.50}}
\put(35.50,10.00){\circle{1.50}}
\put(14.50,10.00){\circle{1.50}}
%%and stars
\put(14.00,6.50){\line(1,1){7.00}}
\put(14.00,13.50){\line(1,-1){7.00}}
\put(29.00,6.50){\line(1,1){7.00}}
\put(29.00,13.50){\line(1,-1){7.00}}
%%%%%%%%%%%%%%%%%%%%%%%%%%%%%%%%%%%%%%%%
%0 horizontal line
\put(-8.00,25.00){\line(1,0){66.00}}
%its vertices
\put(10.00,25.00){\circle*{1.50}}
\put(0.00,25.00){\circle*{1.50}}
\put(40.00,25.00){\circle*{1.50}}
\put(50.00,25.00){\circle*{1.50}}
%%%%%%%%%%%%%%%%%%%%%%%%%%%%%%%%%%%%%
%+1 horizontal line
\put(13.00,40.00){\line(1,0){24.00}}
%its vertices
\put(17.50,40.00){\circle{1.50}}
\put(32.50,40.00){\circle{1.50}}
\put(15.50,38.00){\circle*{1.50}}
\put(15.50,42.00){\circle{1.50}}
\put(19.50,38.00){\circle{1.50}}
\put(19.50,42.00){\circle*{1.50}}
\put(30.50,38.00){\circle{1.50}}
\put(30.50,42.00){\circle*{1.50}}
\put(34.50,38.00){\circle*{1.50}}
\put(34.50,42.00){\circle{1.50}}
\put(35.50,40.00){\circle{1.50}}
\put(14.50,40.00){\circle{1.50}}
%and stars
\put(14.00,36.50){\line(1,1){7.00}}
\put(14.00,43.50){\line(1,-1){7.00}}
\put(29.00,36.50){\line(1,1){7.00}}
\put(29.00,43.50){\line(1,-1){7.00}}
%%%%%%%%%%%%%%%%%%%%%%%%%%%%%%%%%%%%%%
%%+2 horizontal line
\put(15.00,50.00){\line(1,0){20.00}}
%%with vertices
\put(20.00,50.00){\circle{1.50}}
\put(25.00,50.00){\circle*{1.50}}
\put(30.00,50.00){\circle{1.50}}
\put(18.00,48.00){\circle*{1.50}}
\put(18.00,52.00){\circle{1.50}}
\put(22.00,48.00){\circle{1.50}}
\put(22.00,52.00){\circle*{1.50}}
\put(28.00,48.00){\circle{1.50}}
\put(28.00,52.00){\circle*{1.50}}
\put(32.00,48.00){\circle*{1.50}}
\put(32.00,52.00){\circle{1.50}}
\put(16.50,50.00){\circle{1.50}}
\put(33.50,50.00){\circle{1.50}}
%%and stars
\put(46.50,16.50){\line(1,1){7.00}}
\put(46.50,23.50){\line(1,-1){7.00}}
\put(46.50,26.50){\line(1,1){7.00}}
\put(46.50,33.50){\line(1,-1){7.00}}
%%%%%%%%%%%%%%%%%%%%%%%%%%%%%%%%%%%%%
%%%%%%%%%vertical lines %%%%%%%%%%%%%%%%%%%%%%%%%%%%%%%%
%-2 vertical line
\put(0.00,15.00){\line(0,1){20.00}}
%its vertices
\put(0.00,20.00){\circle{1.50}}
\put(0.00,30.00){\circle{1.50}}
\put(-2.00,18.00){\circle{1.50}}
\put(-2.00,22.00){\circle*{1.50}}
\put(2.00,18.00){\circle*{1.50}}
\put(2.00,22.00){\circle{1.50}}
\put(-2.00,28.00){\circle*{1.50}}
\put(-2.00,32.00){\circle{1.50}}
\put(2.00,28.00){\circle{1.50}}
\put(2.00,32.00){\circle*{1.50}}
\put(0.00,33.50){\circle{1.50}}
\put(0.00,16.50){\circle{1.50}}
%and stars
\put(-3.50,16.50){\line(1,1){7.00}}
\put(-3.50,23.50){\line(1,-1){7.00}}
\put(-3.50,26.50){\line(1,1){7.00}}
\put(-3.50,33.50){\line(1,-1){7.00}}
%%%%%%%%%%%%%%%%%%%%%%%%%%%%%%%%%%%%%%%%
%-1 vertical line
\put(10.00,13.00){\line(0,1){24.00}}
%its vertices
\put(10.00,32.50){\circle{1.50}}
\put(10.00,17.50){\circle{1.50}}
\put(8.00,15.50){\circle{1.50}}
\put(8.00,19.50){\circle*{1.50}}
\put(12.00,15.50){\circle*{1.50}}
\put(12.00,19.50){\circle{1.50}}
\put(8.00,30.50){\circle*{1.50}}
\put(8.00,34.50){\circle{1.50}}
\put(12.00,30.50){\circle{1.50}}
\put(12.00,34.50){\circle*{1.50}}
\put(10.00,35.50){\circle{1.50}}
\put(10.00,14.50){\circle{1.50}}
%and stars
\put(6.50,14.00){\line(1,1){7.00}}
\put(6.50,21.00){\line(1,-1){7.00}}
\put(6.50,29.00){\line(1,1){7.00}}
\put(6.50,36.00){\line(1,-1){7.00}}
%%%%%%%%%%%%%%%%%%%%%%%%%%%%%%%%%%%%%%
%0 vertical line
\put(25.00,-8.00){\line(0,1){66.00}}
%its vertices
\put(25.00,10.00){\circle*{1.50}}
\put(25.00,25.00){\circle*{1.50}}
\put(25.00,40.00){\circle*{1.50}}
%%%%%%%%%%%%%%%%%%%%%%%%%%%%%%%%%
%+1 vertical line
\put(40.00,13.00){\line(0,1){24.00}}
%%its vertices
\put(40.00,17.50){\circle{1.50}}
\put(40.00,32.50){\circle{1.50}}
\put(38.00,15.50){\circle*{1.50}}
\put(38.00,19.50){\circle{1.50}}
\put(42.00,15.50){\circle{1.50}}
\put(42.00,19.50){\circle*{1.50}}
\put(38.00,30.50){\circle{1.50}}
\put(38.00,34.50){\circle*{1.50}}
\put(42.00,30.50){\circle*{1.50}}
\put(42.00,34.50){\circle{1.50}}
\put(40.00,35.50){\circle{1.50}}
\put(40.00,14.50){\circle{1.50}}
%%and stars
\put(36.50,14.00){\line(1,1){7.00}}
\put(36.50,21.00){\line(1,-1){7.00}}
\put(36.50,29.00){\line(1,1){7.00}}
\put(36.50,36.00){\line(1,-1){7.00}}
%%%%%%%%%%%%%%%%%%%%%%%%%%%%%%%%%%%%%%%
%%+2 vertical line
\put(50.00,15.00){\line(0,1){20.00}}
%%its vertices
\put(50.00,20.00){\circle{1.50}}
\put(50.00,30.00){\circle{1.50}}
\put(48.00,18.00){\circle*{1.50}}
\put(48.00,22.00){\circle{1.50}}
\put(52.00,18.00){\circle{1.50}}
\put(52.00,22.00){\circle*{1.50}}
\put(48.00,28.00){\circle{1.50}}
\put(48.00,32.00){\circle*{1.50}}
\put(52.00,28.00){\circle*{1.50}}
\put(52.00,32.00){\circle{1.50}}
\put(50.00,16.50){\circle{1.50}}
\put(50.00,33.50){\circle{1.50}}
%%with stars
\put(16.50,46.50){\line(1,1){7.00}}
\put(23.50,46.50){\line(-1,1){7.00}}
\put(26.50,46.50){\line(1,1){7.00}}
\put(33.50,46.50){\line(-1,1){7.00}}
%%root
\put(26.00,26.00){\footnotesize $e$}
\end{picture}
\caption{Matricially free analog of ${\mathbb H}_{4}$}
\end{figure}
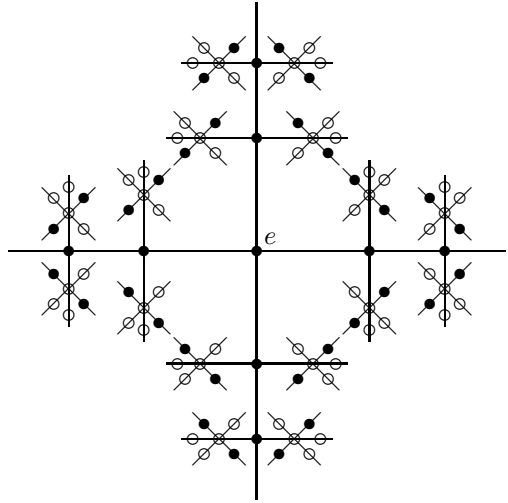
\begin{Example}
{\rm In the above example take an $n$-dimensional square array of 
copies of $F(1)$ and denote their generators by $g_{i,j}$, where $i,j\in[n]$.
For each natural $n$ we form a tree (a subtree of the 
homogenous tree ${\mathbb H}_{2n^{2}}$) which corresponds to 
the `matricially free product of $n^2$ free groups'. 
Suppose the root $e$ corresponds to the `father'.
We distinguish `sons' and `daughters' in each `generation' which 
correspond to the left action of $g_{j,j}$ or $g_{j,j}^{-1}$, 
and $g_{i,j}$ or $g_{i,j}^{-1}$, respectively, where $i\neq j$.
The rules of drawing the tree follow from matricial freeness
and are the following: each `son' has $1$ `son' and $2n-2$ `daughters', whereas
each `daughter' has 2 `sons' and $2n-1$ `daughters'. Therefore, 
`sons' and `daughters' correspond to vertices of valencies $2n$ and $2n+2$, 
respectively. In Figure 1 we draw such a tree for $n=2$ (black
and empty circles are assigned to `sons' and `daughters', respectively).
If we make an additional assumption, for instance that `daughters' cannot have `sons' 
(this fact corresponds to strong matricial freeness, where diagonal generators kill
words beginning with the non-diagonal ones), we recover ${\mathbb H}_{2n}$ of free probability (cf. Remark 3.1).}
\end{Example}
\begin{Example}
{\rm Let
${\mathcal R}(\widehat{\mathcal H})$ be the discrete strongly matricially free Fock space and let
${\mathcal R}(\widehat{\mathcal H})\cong {\mathcal F}(\bigoplus_{j}{\mathbb C}e_{j})$ be 
the natural isomorphism of Remark 3.1, where $\{e_{j}:j\in {\mathbb N}\}$ is an orthonormal basis of some Hilbert space.
If $(w_{i,j})$ is an infinite matrix with non-negative parameters $p$ and $q$ above and below the main diagonal, 
respectively, and $1$'s on the diagonal, then this it is easy to see that the 
$(p,q)$-creation operators studied in [15] can be identified (with the use of this isomorphism) 
with the strongly convergent sums 
$$
A_{i}=\sum_{j}w_{i,j}k_{i,j}
$$ 
where $i\in {\mathbb N}$ and the $(p,q)$-annihilation operators are their adjoints.
A similar approach can be applied to square arrays of arbitrary Hilbert spaces and *-representations, which leads to some notion of `$(p,q)$-independence'. Moreover, it can be carried out for more general matrices $(w_{i,j})$ within the framework of the strong matricial freeness, which generalizes notions of independence of this type.}
\end{Example}

\section{Traces}
In this Section we introduce some real-valued functions on the 
set of non-crossing pair partitions. These functions are obtained by computing traces
of a square real-valued matrix. We will assume later that this matrix
has non-negative entries which represent variances of probability measures on the real line
${\mathcal M}_{{\mathbb R}}$ and we will demonstrate that the functions 
introduced in this Section describe the asymptotics of matricially free 
random variables in central limit theorems.

Let $\mathcal{NC}_{m}$ denote the set of non-crossing partitions of the set $[m]$, 
i.e. if $\pi =\{\pi_1, \pi_2,\ldots , \pi_k\}\in \mathcal{NC}_{m}$, then there are no numbers $i<p<j<q$ such that $i,j\in \pi_{r}$ and $p,q\in \pi_{s}$
for $r\neq s$. The block $\pi_r$ is {\it inner} with respect to $\pi_s$
if $p<i<q$ for any $i\in \pi_r$ and $p,q\in \pi_s$
(then $\pi_s$ is {\it outer} with respect to $\pi_r$).
It is clear that if $\pi_r$ has outer blocks, then there exists a 
unique block among them, say $\pi_s$, which is {\it nearest} to $\pi_r$, 
i.e. if another block, say $\pi_t$,
is outer with respect to $\pi_r$, then we must have $a<p<b$ for 
any $a,b\in \pi_t$ and $p\in \pi_s$.
Then the pair $(\pi_{r},\pi_{s})$ is called the {\it nearest inner-outer pair} of blocks.

Let $\mathcal{NCC}_{m}$ denote the set of {\it non-crossing covered pair partitions} of $[m]$, 
by which we understand the subset of $\mathcal{NC}_{m}$
consisting of those partitions in which $1$ and $m$ 
belong to the same block (if $m=1$, we understand that 
the partition consists of one block).
In terms of diagrams, all blocks of $\pi\in \mathcal{NCC}_{m}$,
where $m>1$, are covered by the block containing $1$ and $m$. 
We denote by $\mathcal{NC}_{m}^{2}$ and $\mathcal{NCC}_{m}^{2}$
the sets of non-crossing pair partitions of $[m]$ and non-crossing covered
pair-partitions of $[m]$, respectively, and we set 
$\mathcal{NC}^{2}=\bigcup_{m=1}^{\infty}\mathcal{NC}_{m}^{2}$ and
$\mathcal{NCC}^{2}=\bigcup_{m=1}^{\infty}\mathcal{NCC}_{m}^{2}$.

It is easy to see that each $\pi\in \mathcal{NC}_{m}$ can be decomposed as 
\begin{equation}
\pi=\pi^{(1)}\cup \pi^{(2)}\cup \ldots\cup \pi^{(p)}
\end{equation}
where $\pi^{(1)}, \ldots , \pi^{(p)}$ are non-crossing covered partitions of 
subintervals $I_{1},I_{2}, \ldots , I_{p}$ of $[m]$ whose union
gives $[m]$. By a partition of a set $I$ consisting of $r$ elements 
we understand the corresponding partition of $[r]$.

On the other hand, each $\pi\in \mathcal{NCC}_{m}$ can be decomposed as
\begin{equation}
\pi=\pi^{(0)}\cup \pi^{(1)}\cup \ldots \cup \pi^{(r)}
\end{equation}
where $\pi^{(0)}$ is the block containing $1$ and $m$ 
and $\pi^{(1)}, \pi^{(2)}, \ldots , \pi^{(r)}$ 
are non-crossing covered partitions of 
subintervals $I_{1},I_{2}, \ldots , I_{r}$ of the set 
$\{2, \ldots , m-1\}$.

Consider now a square real-valued matrix $V=(v_{i,j})\in M_{n}({\mathbb R})$, 
where $n\in {\mathbb N}\cup \{\infty\}$. 
The usual trace and the normalized trace will be denoted
$$
{\rm Tr}(V)=\sum_{j=1}^{n}v_{j,j} \;\;\; {\rm and} \;\;\; 
{\rm tr}(V)=\frac{1}{n}\sum_{j=1}^{n}v_{j,j},
$$
respectively.
For finite $n$, we define the `diagonalization mapping'
$$
\tau:\; M_{n}({\mathbb R})\rightarrow D_{n}({\mathbb R}), \;\;\;
\tau(V)={\rm diag}(\sum_{j}v_{j,1}, \ldots , \sum_{j}v_{j,n})
$$
where $D_{n}({\mathbb R})$ is the set of square diagonal real-valued matrices of dimension $n$
(by abuse of notation, the same symbol $\tau$ is used for all $n$).
In other words, $\tau$ computes the sum of all elements of $V$ in each column separately and puts 
this value on the diagonal.

Using the above trace operations, we shall define two real-valued functions on the set $\mathcal{NC}^{2}$, denoted $v$ and $v_{0}$,  associated with given 
$V\in M_{n}({\mathbb R})$. Although there is a close similarity between these functions
when restricted to $\mathcal{NCC}^{2}$ (in particular, they are defined as traces of certain 
matrix-valued quasi-mulitplicative functions), note that they are extended to $\mathcal{NC}^{2}$ in two different ways.
\begin{Definition}
{\rm For a given matrix $V\in M_{n}({\mathbb R})$, we define a mapping
from $\mathcal{NC}^{2}$ to $D_{n}({\mathbb R})$ by assigning to 
each $\pi\in \mathcal{NC}^{2}$ the matrix $V(\pi)$ by the following recursion:
\begin{enumerate}
\item
if $\pi$ consists of one block, we set $V(\pi)=\tau(V)$,
\item
if $\pi\in \mathcal{NCC}^{2}$ consists of more than one block, then
\begin{equation}
V(\pi)=\tau(V(\pi^{(1)})\ldots V(\pi^{(r)})V)
\end{equation}
according to the decomposition (5.2),
\item
if $\pi\in \mathcal{NC}^{2}$, then
\begin{equation}
V(\pi)=V(\pi^{(1)})V(\pi^{(2)})\ldots V(\pi^{(p)})
\end{equation}
according to the decomposition (5.1).
\end{enumerate}
Let $v:\mathcal{NC}^{2}\rightarrow {\mathbb R}$ be the function defined by
$v(\pi)={\rm tr}(V(\pi))$.}
\end{Definition}
\begin{Example}
{\rm 
For some $V\in M_{n}({\mathbb R})$, consider three partitions given in Fig.2.
Color each block by a number from the set $[n]$.
Computation of the corresponding values of the function $v$ gives
\begin{eqnarray*}
v(\pi)&=&{\rm tr} (\tau(\tau(V)V))=\frac{1}{n}\sum_{i,j,k}v_{i,j}v_{j,k}\\
v(\eta)&=&{\rm tr}
\left(\tau\left(\tau(V)V\right)\tau(V)\right)=\frac{1}{n}\sum_{i,j,k,l}v_{i,j}v_{j,l}v_{k,l}\\
v(\zeta)&=&{\rm tr} (\tau(\tau(V)\tau(V)V))=\frac{1}{n}\sum_{i,j,k,l}v_{i,k}v_{j,k}v_{k,l}
\end{eqnarray*}
One can see that to each nearest inner-outer pair of blocks 
$(\pi_r,\pi_s)$ we assign the matrix element $v_{p,q}$, where $p$ and $q$ are the colors of
$\pi_r$ and $\pi_s$ respectively.
Moreover, if a block does not have any outer blocks and is colored by $q$, 
then we assign to it the matrix element $v_{q,t}$, where 
$t$ is assumed to be the same for all such blocks
(one can imagine that we have an additional `conditional block' colored by $t$ 
which covers all other blocks). At the end we sum over all colorings.}
\end{Example}
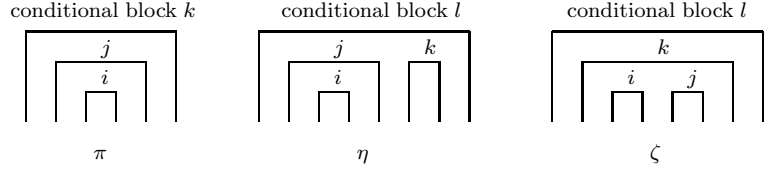
\begin{figure}
\unitlength=1mm
\special{em:linewidth 0.4pt}
\linethickness{0.4pt}
\begin{picture}(120.00,30.00)(-15.00,5.00)
%%%%%%%%%% diagram 1 %%%%%%%%%%%%%%%%%%%%
\put(-5.00,10.00){\line(0,1){12.00}}
\put(-1.00,10.00){\line(0,1){8.00}}
\put(3.00,10.00){\line(0,1){4.00}}
\put(7.00,10.00){\line(0,1){4.00}}
\put(11.00,10.00){\line(0,1){8.00}}
\put(15.00,10.00){\line(0,1){12.00}}

\put(-7.00,24.00){$\scriptstyle{{\rm conditional\;block}\;k}$}
\put(5.00,19.00){$\scriptstyle{j}$}
\put(5.00,15.00){$\scriptstyle{i}$}

\put(-5.00,22.00){\line(1,0){20.00}}
\put(-1.00,18.00){\line(1,0){12.00}}
\put(3.00,14.00){\line(1,0){4.00}}
\put(4.00,5.00){$\scriptstyle{\pi}$}
%%%%%%%%%%%%diagram 2%%%%%%%%%%%%%%%%%%%%%%%
\put(26.00,10.00){\line(0,1){12.00}}
\put(54.00,10.00){\line(0,1){12.00}}
\put(30.00,10.00){\line(0,1){8.00}}
\put(34.00,10.00){\line(0,1){4.00}}
\put(38.00,10.00){\line(0,1){4.00}}
\put(42.00,10.00){\line(0,1){8.00}}
\put(46.00,10.00){\line(0,1){8.00}}
\put(50.00,10.00){\line(0,1){8.00}}

\put(29.00,24.00){$\scriptstyle{{\rm conditional \;block}\;l}$}
\put(36.00,19.00){$\scriptstyle{j}$}
\put(36.00,15.00){$\scriptstyle{i}$}
\put(48.00,19.00){$\scriptstyle{k}$}

\put(26.00,22.00){\line(1,0){28.00}}
\put(30.00,18.00){\line(1,0){12.00}}
\put(34.00,14.00){\line(1,0){4.00}}
\put(46.00,18.00){\line(1,0){4.00}}
\put(39.00,5.00){$\scriptstyle{\eta}$}
%%%%%%%%%%%diagram 3%%%%%%%%%%%%%%%%%%%%%%%%%%%%%
\put(65.00,10.00){\line(0,1){12.00}}
\put(69.00,10.00){\line(0,1){8.00}}
\put(73.00,10.00){\line(0,1){4.00}}
\put(77.00,10.00){\line(0,1){4.00}}
\put(81.00,10.00){\line(0,1){4.00}}
\put(85.00,10.00){\line(0,1){4.00}}
\put(89.00,10.00){\line(0,1){8.00}}
%\put(93.00,10.00){\line(0,1){8.00}}
%\put(97.00,10.00){\line(0,1){8.00}}
\put(93.00,10.00){\line(0,1){12.00}}

\put(75.00,15.00){$\scriptstyle{i}$}
\put(83.00,15.00){$\scriptstyle{j}$}
\put(79.00,19.00){$\scriptstyle{k}$}
%\put(95.00,19.00){$\scriptstyle{l}$}
\put(67.00,24.00){$\scriptstyle{{\rm conditional\;block}\;l}$}

\put(65.00,22.00){\line(1,0){28.00}}
\put(69.00,18.00){\line(1,0){20.00}}
%\put(93.00,18.00){\line(1,0){4.00}}
\put(73.00,14.00){\line(1,0){4.00}}
\put(81.00,14.00){\line(1,0){4.00}}
\put(78.00,5.00){$\scriptstyle{\zeta}$}
%%%%%%%%%%%%%%%%%%%%%%%%%%%%%%%%%%%%%%%%%%%

\end{picture}
\caption{Colored partitions with conditional blocks.}
\end{figure}
The function $v_{0}$ is defined on the set $\mathcal{NCC}^{2}$
in a very similar manner, except that we replace the right multiplier $V$ in (5.3) 
by its main diagonal
$$
V_{0}:={\rm diag}(v_{1,1},\ldots, v_{n,n})
$$
which corresponds to changing only the contribution of the covering block.
Then we extend $v_{0}$ to all of $\mathcal{NC}^{2}$ by multiplicativity
of $v_{0}$.
\begin{Definition}
{\rm For a given matrix $V\in M_{n}({\mathbb R})$, we define
a mapping from $\mathcal{NCC}^{2}$ to $D_{n}({\mathbb R})$
by assigning to each $\pi\in\mathcal{NCC}^{2}$ the matrix $V_{0}(\pi)$
by
 the following recursion:}
\begin{enumerate}
\item
{\rm if $\pi$ consists of one block, we set
$V_{0}(\pi)=V_{0}$}
\item
{\rm if $\pi\in \mathcal{NCC}^{2}$ consists of more than one block, then}
\begin{equation}
V_{0}(\pi)=\tau(V(\pi^{(1)})\ldots V(\pi^{(r)})V_{0})
\end{equation}
{\rm according to the decomposition (5.2).}
\end{enumerate}
{\rm Let $v_{0}:\mathcal{NC}^{2}\rightarrow {\mathbb R}$
be the function defined by $v_{0}(\pi)={\rm Tr}(V_{0}(\pi))$ for $\pi\in \mathcal{NCC}^{2}$,
extended to $\mathcal{NC}^{2}$ by mulitplicativity
$v_{0}(\pi)=v_{0}(\pi^{(1)})\ldots v_{0}(\pi^{(p)})$ 
according to the decomposition (5.1).}
\end{Definition}
\begin{Example}
{\rm 
Let us compute the values of $v_{0}$ corresponding to the partitions in Fig.2.
We have
\begin{eqnarray*}
v_{0}(\pi)&=&{\rm Tr}(\tau(V)V_{0})=\sum_{i,j}v_{i,j}v_{j,j}\\
v_{0}(\eta)&=&{\rm Tr}(\tau(V)V_{0}){\rm Tr}(V_{0})=
\sum_{i,j,k}v_{i,j}v_{j,j}v_{k,k}\\
v_{0}(\zeta)&=&{\rm Tr}
(\tau(\tau(V)\tau(V)V_{0}))=\sum_{i,j,k}v_{i,k}v_{j,k}v_{k,k}
\end{eqnarray*}
Note that the main difference (apart from normalization) between $v_{0}(\pi)$ and $v(\pi)$ 
concerns the blocks which do not have outer blocks. Here, if such a block 
is colored by $q$, we assign to it the matrix element $v_{q,q}$ and we do not use
`conditional blocks'.}
\end{Example}
Below we shall prove a lemma which gives explicit combinatorial formulas
for $v(\pi)$ and $v_{0}(\pi)$.
If $\pi=\{\pi_1, \ldots , \pi_k\}\in \mathcal{NC}_{2k}^{2}$, we color 
the blocks of $\pi$ and the conditional block $\pi_{0}$ by the set $[n]$. 
The {\it coloring} of the partition $\pi\cup \{\pi_{0}\}$ can then be identified with 
the mapping $f:[k]\cup \{0\}\rightarrow [n]$, where $f(i)$ is interpreted as the color of 
$\pi_i$. The pair $(\pi,f)$ can then be interpreted as 
the colored partition associated with $\pi$
and the mapping $f$. The set of all colorings of $\pi$ by the set $[n]$ will be denoted $F_n(\pi)$.

For a given matrix $V=(v_{i,j})\in M_{n}({\mathbb R})$ and given 
colored partition $(\pi,f)$, where $\pi$ is a non-crossing pair partition consisting of $k$ blocks
and $f$ is the coloring of $\pi$, we can now define (slightly abusing notation) 
a natural mapping
$$
v_{0}:\{(\pi_1,f), \ldots , (\pi_k,f)\}\rightarrow {\mathbb R},
$$
by the following rule:
$$
v_0(\pi_p,f)=v_{i,j}\;\;{\rm if}\;\;f(p)=i\;\;{\rm and}\;\;f(\sigma(p))=j\\
$$
where $\sigma(p)$ is the index of the nearest outer block of $\pi_p$ if $\pi_p$ has an outer block,
and otherwise $\sigma(p)=p$. In other words, to each block
we assign the matrix element whose first index is the color of that block, and the second index --
the color of its nearest outer block.

The multiplicative extension of this function to the class of non-crossing colored partitions
gives
$$
v_{0}(\pi,f)=v_0(\pi_1,f)v_0(\pi_2,f)\ldots v_0(\pi_k,f).
$$
In a similar way we define functions
$$
v_{j}(\pi,f)=v_j(\pi_1,f)v_j(\pi_2,f)\ldots v_j(\pi_k,f)
$$
for any $j\in [n]$, where $v_j(\pi_p,f)=v_{i,j}$ if $f(p)=i$ and $f(\sigma_j(p))=j$
where $\sigma_j(p)$ is the index of the nearest outer block of $\pi_p$ if $\pi_p$ has an outer block,
and otherwise we put $\sigma_j(p)=j$. In the latter case, $j$ is interpreted as the color of the conditional block. 
\begin{Lemma}
For any $V\in M_{n}({\mathbb R})$, where $n\in {\mathbb N}$, and $\pi\in \mathcal{NC}^{2}$,
it holds that
$$
v_{0}(\pi)=\sum_{f\in F_n(\pi)}v_{0}(\pi,f)
\;\;\;
and
\;\;\; 
v(\pi)=\frac{1}{n}\sum_{f\in F_n(\pi)}\sum_{j \in [n]}v_{j}(\pi,f)
$$
where the summation over $j$ corresponds to all colorings of the conditional block.
\end{Lemma}
{\it Proof.}
We provide an induction proof for the function $v$ (the proof for $v_{0}$ is similar). 
The main induction step will be carried out on the level of 
the (diagonal) matrices $V(\pi)$. We claim that its diagonal entries are of the form
$$
(V(\pi))_{q,q}= \sum_{f\in F_n(\pi)}v_q(\pi,f)
$$
for any $\pi \in \mathcal{NCC}^{2}$ and $q\in [n]$. In view of (5.4), 
the required formula for $v(\pi)$ is then a straightforward consequence of
the claim. Of course, if $\pi$ consists of one block, then 
$$
(V(\pi))_{q,q}=\sum_{j}v_{j,q}
$$
and thus our assertion easily follows.
Assume now that $\pi$ has $k>2$ blocks and suppose the assertion holds for 
non-crossing covered partitions which have less than $k$ blocks.
Since $\pi$ has a decomposition of type (5.2), the assertion holds for
$\pi^{(1)}, \pi^{(2)}, \ldots , \pi^{(r)}$ in this decomposition.
We know that the matrix assigned to $\pi$ has the form (5.3)
and therefore the product of (diagonal) matrices
corresponding to these subpartitions has (diagonal) matrix elements of the form
$$
(V(\pi))_{q,q}=(\tau(WV))_{q,q}=\sum_{j}W_{j,j}v_{j,q}
$$
where $q$ is the color of the conditional block of $\pi$ and 
$$
W_{j,j}
=
\sum_{f_1\in F_n(\pi^{(1)})} v_j(\pi^{(1)},f_1)\;\;\ldots 
\sum_{f_r\in F_n(\pi^{(r)})} v_j(\pi^{(r)},f_r)
$$
by the inductive assumption.
Now, since the blocks of $\pi^{(1)}, \ldots , \pi^{(r)}$ are 
colored independently, we have
$$
v_j(\pi^{(1)},f_1)\ldots v_j(\pi^{(r)},f_r)v_{j,q}=v_{q}(\pi,f)
$$
for a uniquely determined coloring $f$ of the blocks of $\pi$ in which $j$ can be interpreted
as the color of $\pi^{(0)}$ since $\pi^{(0)}$ covers $\pi^{(1)}, \ldots , \pi^{(r)}$
and $q$ can be viewed as the color of the conditional block of $\pi$. 
This and the above formula for $W_{j,j}$ gives the desired formula
$$
V(\pi)_{q,q}=\sum_{f\in F_n(\pi)}v_{q}(\pi,f),
$$
which proves our claim and thus completes the proof of the theorem.
\hfill $\blacksquare$\\
\indent{\par}
Under suitable assumptions on $V$, there is a simple connection between functions $v$ and $v_{0}$ 
if $V_{0}=aI_{n}/n$, where $I_{n}$ is a unit matrix and $a$ is a positive number.
We shall express this relation in terms of the corresponding
formal Laurent series, which turn out to be Cauchy transforms
of (compactly supported) probability measures on the real line 
associated with appropriately constructed random 
variables,
\begin{Proposition}
Let $V\in M_{n}({\mathbb R})$ be such that $V_{0}=aI_{n}/n$, where $a> 0$.
If $G(z)=\sum_{k=0}^{\infty}a_{2k}z^{-2k-2}$ and $G_{0}(z)=\sum_{k=0}^{\infty}b_{2k}z^{-2k-2}$
are formal Laurent series, where
$$
a_{2k}=\sum_{\pi\in \mathcal{NC}^{2}_{2k}}v(\pi)\;\;\;{\rm and}\;\;\;
b_{2k}=\sum_{\pi\in \mathcal{NC}^{2}_{2k}}v_{0}(\pi)
$$
and both $v(\pi)$ and $v_{0}(\pi)$ are associated with $V$, then $G_{0}(z)=1/(z-aG(z))$.
\end{Proposition}
{\it Proof.}
Observe that in the case when $v_{j,j}=a/n$ for all $j \in [n]$, we have
$$
v(\pi)=\frac{v_{0}(\pi')}{a}
$$
for any $\pi\in \mathcal{NC}_{m}^{2}$, where the partition 
$\pi'\in \mathcal{NCC}_{m+2}^{2}$ is obtained from $\pi$ 
by adding to $\pi$ the block that covers all blocks of $\pi$, say $\{0,m+1\}$.
This leads to 
\begin{eqnarray*}
A(z):&=&\sum_{m=0}^{\infty}a_{2m}z^{2m}\;=\;1+\sum_{m=1}^{\infty}\sum_{\pi\in \mathcal{NC}_{2m}^{2}}v(\pi)\;
z^{2m}\\
&=& 1+ \frac{1}{a}\sum_{m=1}^{\infty}\sum_{\pi\in \mathcal{NCC}_{2m+2}^{2}}v_{0}(\pi)z^{2m}
\;=\; \frac{C(z)-1}{az^2}
\end{eqnarray*}
where $C(z)=\sum_{m=0}^{\infty}c_{2m}z^{2m}$ and 
$c_{2m}=\sum_{\pi\in \mathcal{NCC}_{2m}^{2}}v_{0}(\pi)$. 
Now, using the multiplicativity of $v_{0}$,
we obtain
$$
B(z):=\sum_{m=0}^{\infty}b_{2m}z^{2m}=1+\sum_{m=1}^{\infty}(C(z)-1)^{m}=\frac{1}{2-C(z)}
$$
which leads to 
$$
A(z)=\frac{1}{az^{2}}\left(1-\frac{1}{B(z)}\right)\;\;\;{\rm and}\;\;\;
G_{0}(z)=\frac{1}{z-aG(z)}
$$
where 
$$
G(z)=\frac{1}{z}A\left(\frac{1}{z}\right)\;\;\;{\rm and}\;\;\;
G_{0}(z)=\frac{1}{z}B\left(\frac{1}{z}\right)
$$
which completes the proof.
\hfill $\blacksquare$

\section{Random pseudomatrices}
In this Section we will study the asymptotic behavior of random pseudomatrices for two sequences of states:
the sequence of distinguished states, with respect to which our variables will be matricially free,
and the sequence of traces which are normalized sums of conditions in the definition of matricial freeness.
Under suitable assumptions, we will later obtain two types of central limit theorems:
the standard central limit theorem as well as
the tracial central limit theorem for matricially free random variables, the latter
being related to random matrix models. 

Let $({\mathcal A}_{n})_{n\in {\mathbb N}}$ be a sequence of unital *-algebras. 
For each $n$, let $(X_{i,j}(n))_{1 \leq i,j\leq n}$ 
be an array of self-adjoint random variables in ${\mathcal A}_{n}$ and let
$(\phi_{i,j}(n))$ be an array of associated states on ${\mathcal A}_{n}$.
We will say that the moments of the $X_{i,j}(n)$ 
are {\it uniformly bounded} with respect to $(\phi_{i,j}(n))$ if, for all natural $m$, 
there exists $M_{m}\geq 0$ such that 
\begin{equation}
|\phi_{i,j}(n)(X_{i,j}^{m}(n))|\leq \frac{M_{m}}{n^{m/2}}
\end{equation}
for all $n\in {\mathbb N}$ and $i,j\in [n]$.

Assume further that for each natural $n$ the array $(X_{i,j}(n))$ is matricially free
with respect to the array $(\phi_{i,j}(n))$ defined by a distinguished state 
$\phi_n$ and the conditions $\psi_{n,j}$, where $1\leq j \leq n$. 
We are going to study the asymptotic behavior of random pseudomatrices 
\begin{equation}
S(n)=\sum_{i,j=1}^{n}X_{i,j}(n)
\end{equation}
with respect to two types of states:
\begin{enumerate}
\item
distinguished states $\phi_{n}$,
\item
convex linear combinations of conditions
\begin{equation}
\psi_{n}:=\frac{1}{n}\sum_{j=1}^{n}\psi_{n,j}
\end{equation} 
\end{enumerate}
which play the role of traces in the context of random pseudomatrices.

This will lead to two types of central limit theorems for matricially free random variables: 
the `standard CLT' reminding the CLT for free random variables [20] and the `tracial CLT' reminding
the limit theorem for random matrices [22]. Thanks to the `matricial property' and the
`diagonal subordination property' of the matricially free product of states, the normalization 
of square root type works in both cases since the number of summands in $S(n)$ 
which give a non-zero contribution to the limits is in both cases of order $n$.

Of course, the uniform boundedness assumption is satisfied if we take variables of type $X_{i,j}(n)=X_{i,j}/\sqrt{n}$, where $(X_{i,j})$ is an infinite array of random variables
whose distributions in the states $\phi_{i,j}(n)$, respectively, are identical and do not depend on $n$.
Although it is convenient to think of the $X_{i,j}(n)$ as if they were of this form, 
we want to study similar variables, whose variances are of order $1/n$ and stay the same within blocks whose sizes become infinite as $n\rightarrow \infty$.

If the tuple $((i_1,j_1), \ldots, (i_m,j_m))\in (I\times I)^{m}$, where
$I$ is an index set, defines a partition $\pi=\{\pi_1,\ldots, \pi_k\}$ 
of the set $[m]$, i.e. $(i_p,j_p)=(i_q,j_q)$
if and only if there exists $r$ such that $p,q\in \pi_{r}$, we will write
$$
P((i_1,j_1),\ldots, (i_m,j_m))=\pi \;.
$$ 
Of course, if $\pi$ is a non-crossing pair partition, then each $\pi_{r}$ is a
two-element set. 
We will also adopt the convention that if $m$ is odd, then
$\mathcal{NC}_{m}^{2}=\emptyset$ and the summation over $\pi\in \mathcal{NC}_{m}^{2}$
gives zero. This allows us to state results for moments of the $S(n)$ of all orders 
without distinguishing even and odd moments.

\begin{Lemma}
Let $(\mathcal{A}_{n},\phi_{n})$ be a sequence noncommutative probability spaces, each with 
an array $(\phi_{i,j}(n))$ defined as above and such that 
\begin{enumerate}
\item
for fixed $n\in {\mathbb N}$, the array $(X_{i,j}(n))$ is 
matricially free with respect to $(\phi_{i,j}(n))$,
\item 
$\phi_{i,j}(n)(X_{i,j}(n))=0$ and $\phi_{i,j}(n)(X^{2}_{i,j}(n))=v_{i,j}(n)$ 
for any $i,j$ and $n\in {\mathbb N}$,
\item
the moments of the $X_{i,j}(n)$ are uniformly bounded with respect to $(\phi_{i,j}(n))$.
\end{enumerate}
Then
$$
\psi_{n}(S^{m}(n))=
\sum_{\pi\in {\mathcal NC}_{m}^{2}}
v(\pi,n)+ O\left(\frac{1}{\sqrt{n}}\right)
$$
where $v(\pi,n)={\rm tr}(V(\pi,n))$ is given by Definition 5.1 and
corresponds to the variance matrix $V(n)=(v_{i,j}(n))$. 
\end{Lemma}
{\it Proof.}
We have
\begin{eqnarray*}
\psi_{n}(S^{m}(n))&=&
\sum_{i_1,j_1, \ldots , i_m,j_m}
\psi_{n}(X_{i_1,j_1}(n)\ldots X_{i_m,j_m}(n))\\
&=&
\sum_{\pi\in {\mathcal P}_{m}}
\sum_{\stackrel {i_1,j_1, \ldots , i_m,j_m}
{\scriptscriptstyle P((i_1,j_1), \ldots , (i_m,j_m))=\pi}}
\psi_{n}(X_{i_1,j_1}(n)\ldots X_{i_m,j_m}(n))
\end{eqnarray*}
where $\mathcal{P}_{m}$ denotes the set of all partitions of $[m]$.
Standard arguments allow us to conclude that for large $n$ only non-crossing pair partitions
give relevant contributions, namely
\begin{enumerate}
\item
the moments of
the $X_{i,j}(n)$ satisfy the singleton condition with respect to the $\psi_{n,k}$
and thus with respect to $\psi_{n}$, namely
$$
\psi_{n}(X_{i_1,j_1}(n)\ldots X_{i_m,j_m}(n))=0
$$
if $\pi=P((i_1,j_1), \ldots , (i_m,j_m))$ contains a singleton,
\item
if $\pi$ has no singletons and contains blocks consisting of more than two elements, then 
$$
\sum_{\stackrel {i_1,j_1, \ldots , i_m,j_m}
{\scriptscriptstyle P((i_1,j_1), \ldots , (i_m,j_m))=\pi}}
\psi_{n}(X_{i_1,j_1}(n)\ldots X_{i_m,j_m}(n))=
O\left(\frac{1}{\sqrt{n}}\right)
$$
\item
if $m$ is even and $\pi$ is a crossing pair 
partition of $[m]$, then the contribution from the corresponding mixed moments
of the $X_{i,j}$ is zero.
\end{enumerate}
The above arguments are similar to those that would hold if 
random variables $X_{i,j}(n)$, where $(i,j)\in I\times I$, were free.
(1) follows from the definition of the matricially free product of states 
and the mean zero assumption. (2) holds since moments are uniformly bounded and the usual combinatorial argument (there are fewer than $m/2$ independent indices to sum over) can be used.
In particular, this and (1) imply that if $m$ is odd, then the contribution to the limit is of order
$O(1/\sqrt{n})$. (3) is a consequence of the mean zero assumption and the 
definition of the matricial freeness (the mixed moments in our case
vanish always when the corresponding moments vanish in the free case).

Therefore,
$$
\psi_{n}(S^{m}(n))=
\sum_{\pi\in {\mathcal NC}_{m}^{2}}
\sum_{\stackrel {i_1,j_1, \ldots , i_m,j_m}
{\scriptscriptstyle P((i_1,j_1), \ldots , (i_m,j_m))=\pi}}
\psi_{n}(X_{i_1,j_1}(n)\ldots X_{i_m,j_m}(n))
+O\left(\frac{1}{\sqrt{n}}\right).
$$
Now, suppose that $m$ is even, $\pi\in \mathcal{NC}_{m}^{2}$ and
the sequence of pairs $((i_1,j_1), \ldots (i_m,j_m))$ is 
`compatible with the matricial multiplication', i.e. such that
if $(i_k,j_k)$ and $(i_r,j_r)$ label an inner-outer pair of blocks,
then $j_k=i_r$. If $\pi_j=\{r,r+1\}$ is a block which has no inner blocks,
then we can `pull out' the variance corresponding to that block, namely:
$$
\psi_{n}(X_{i_1,j_1}(n)\ldots X_{i_m,j_m}(n))
$$

$$
=v_{i_r,j_{r}}
\psi_{n}(X_{i_1,j_1}(n)\ldots X_{i_{r-1},j_{r-1}}(n)X_{i_{r+2},j_{r+2}}(n)
\ldots X_{i_{m},j_m}(n))
$$
(mean zero assumption is used again). 
Continuing this procedure with other blocks which have no inner blocks, 
and summing over indices $i_1,j_1, \ldots , i_m,j_m$, for which it holds that
$P((i_1,j_1), \ldots , (i_m,j_m))=\pi$, we arrive at
$$
\frac{1}{n}
\sum_{k_0,k_1\ldots , k_m}v_{k_1,k_{\sigma(1)}}(n)\ldots v_{k_{m},k_{\sigma(m)}}(n)
+O\left(\frac{1}{\sqrt{n}}\right)
$$
where $\sigma(r)=0$ if $\pi_r$ has no outer blocks ($k_{0}$ labels the conditional block)
and $\sigma(r)=j$ if the nearest outer block of $\pi_r$ is labelled by $j$. 
Let us observe that we included in the above sum all possible labellings of 
the blocks of $\pi$. This is done for convenience since 
it enables us to express the final result in terms of $v(\pi)$. More explicitly, 
we allow $k_0,k_1, \ldots , k_m$ to assume arbitrary values from the set $[m]$
(in particular, they can all be equal), which produces certain terms which cannot be 
obtained from the summation over all $((i_1,j_1), \ldots , (i_m,j_m))$ which define $\pi$.
For example, no nearest inner-outer pair of blocks can contribute $v_{j,j}v_{j,j}$, which appears
in the above sum. However, all such terms are of order $1/\sqrt{n}$ due to insufficient 
number of different summation indices (there are fewer than $m/2$ independent indices) 
and therefore they can be included in the sum without changing the asymptotics.
Using Lemma 5.1, we obtain our assertion.\hfill $\blacksquare$
\begin{Lemma}
Under the assumptions of Lemma 6.1 it holds that
$$
\phi_{n}(S^{m}(n))=
\sum_{\pi\in {\mathcal NC}_{m}^{2}}
v_{0}(\pi,n) + O\left(\frac{1}{\sqrt{n}}\right)
$$
where $v_{0}(\pi,n)={\rm Tr}(V_{0}(\pi,n))$ 
is given by Definition 5.2 and corresponds to the variance matrix $V(n)=(v_{i,j}(n))$.
\end{Lemma}
{\it Proof.}
The proof is similar to that of Lemma 6.1.
\hfill $\blacksquare$\\
\indent{\par}
In order to ensure existence of the limits  of $v(\pi,n)$
and $v_{0}(\pi,n)$ as $n \rightarrow \infty$, we need to  make some assumptions
on the sequence of matrices $(V(n))_{n\in {\mathbb N}}$. For that purpose, introduce a matrix, called
the {\it dimension matrix}
$$
D={\rm diag}(d_1,d_2, \ldots ,d_r)\in D_{r}({\mathbb R})\;\;\;{\rm with}\;\;\;
{\rm Tr}(D)=1,
$$
where $d_{1},d_{2}, \ldots , d_{r}$ are positive real numbers.
For given natural $n$, we associate with $D$ the partition
$[n]=N_{1}\cup N_{2}\cup \ldots \cup N_{r}$, where 
$$
N_1=[1, n_1],N_2=[1+n_1,n_1+n_2],\ldots , N_{r}=[1+n_1+\ldots +n_{r-1}, n]
$$
and
$$
n_k=E\left(\sum_{i=1}^{k}d_{i}n\right)-E\left(\sum_{i=1}^{k-1}d_{i}n\right)
$$
for each $k \in [r]$, with $E(x)$ denoting the largest 
integer smaller or equal to $x$. Of course, the $N_{k}$ are 
disjoint and naturally ordered intervals containing 
$n_{k}$ natural numbers, respectively. Note that in the limit 
$n\rightarrow \infty$ we obtain $n_{k}/n\rightarrow d_{k}$ 
for each $k$.

Assume now that each matrix $V(n)$ has the block-form
\begin{equation}
V(n)=\left(
\begin{array}{ccc}
V_{1,1}(n)  &\ldots & V_{1,r}(n)\\
\vdots & \ddots & \vdots \\
V_{r,1}(n) &\ldots & V_{r,r}(n)
\end{array}
\right)
\end{equation}
where each block $V_{i,j}(n)$ consists of
the same number $u_{i,j}/n$. 
In other words, each $V(n)$ is obtained from a real-valued matrix 
$$
U=(u_{i,j})\in M_{r}({\mathbb R})
$$
by repeating $n_{i}\times n_{j}$ times each entry $u_{i,j}$ 
at all entries of block  $V_{i,j}(n)$ and dividing it by $n$.
Unless stated otherwise, 
we assume that $u_{j,j}>0$ for any $j$ and that $u_{i,j}\geq 0$ for any $i\neq j$.
The whole sequence $(V(n))_{n\in {\mathbb N}}$ 
is built in such a way that proportions between sizes of these blocks are similar for all $n$
and expressed in terms of the dimension matrix $D$ (asymptotically,
these proportions are given by the proportions between numbers $d_{j}$).

Assuming that the variance matrices $V(n)$ are of the above block form, 
we can now state the standard and tracial central limit theorems, 
with limit distributions described in terms of traces of Section 5. 

\begin{Theorem}
Under the assumptions of Lemma 6.1, if $\,V(n)$ is of the block form (6.4) 
for each $n\in {\mathbb N}$, then 
\begin{equation}
\lim_{n \rightarrow \infty}
\psi_{n}(S^{m}(n))=
\sum_{\pi\in {\mathcal NC}_{m}^{2}}b(\pi),
\end{equation}
for any $m\in {\mathbb N}$, where $b(\pi)={\rm Tr}(B(\pi)D)$ and $B(\pi)$ is the
diagonal matrix of Definition 5.1 correspon\-ding to $\pi$ and the matrix $B=DU$.
\end{Theorem}
{\it Proof.}
Clearly, if $m$ is odd, we get zeros on both sides of the above formula (we use our convention
that in this case $\mathcal{NC}_{m}^{2}=\emptyset$). 
The proof for $m=2k$, where $k\in {\mathbb N}$, is based on Lemmas 5.1 and 6.1. 
If, in the combinatorial expression for $v(\pi,n)$,
we substitute for the matrix elements of $V(n)$ the assumed block form,
then, using the partition of the set of colors $[n]=N_1\cup N_2\cup \ldots \cup N_r$,
we can perform summations over the colorings which 
belong to each interval $N_j$ separately. 
Thus, the contributions of various $[n]$-colorings 
of $\pi$ to the limit laws reduce to those corresponding to 
$[r]$-colorings and are described in terms of numbers 
$u_j(\pi_i,f)$, where $i\in [k]$, $j\in [r]$ and $f\in F_{r}(\pi)$ 
(the number $j$ is the color of the conditional block).
We have
\begin{eqnarray*}
\lim_{n \rightarrow \infty}v(\pi,n)&=&\lim_{n\rightarrow \infty}
\left(\frac{1}{n^{k+1}}
\sum_{j\in [r]}n_j\sum_{f\in F_{r}(\pi)}
n_{f(1)}u_j(\pi_1,f)\ldots n_{f(k)}u_j(\pi_k,f)\right)\\
&=&
\sum_{j\in [r]}d_j\sum_{f\in F_{r}(\pi)}
d_{f(1)}u_j(\pi_1,f)\ldots d_{f(k)}u_j(\pi_k,f)\;\;
=\;
{\rm Tr}(B(\pi)D)
\end{eqnarray*}
where $\pi\rightarrow B(\pi)$ is the matrix-valued function which corresponds to 
the matrix $B=DU$ in accordance with Definition 5.1.
In terms of matrix multiplication, the expression on the right hand side
is obtained from that of Lemma 5.1 corresponding to matrix $U$ 
by multiplying $U$ from the left by the dimension 
matrix $D$ and multiplying the whole product of matrices from the right 
by $D$. This proves our assertion.
\hfill $\blacksquare$\\
\begin{Theorem}
Under the assumptions of Lemma 6.2, if $\,V(n)$ is of the block form (6.4) for each 
$n\in {\mathbb N}$, then 
\begin{equation}
\lim_{n \rightarrow \infty}
\phi_{n}(S^{m}(n))=
\sum_{\pi\in {\mathcal NC}_{m}^{2}}
b_{0}(\pi),
\end{equation}
for any $m\in {\mathbb N}$,
where $\pi\rightarrow b_{0}(\pi)$ is the real-valued function of Definition 5.2 
corresponding to the matrix $B=DU$. 
\end{Theorem}
{\it Proof.}
The proof is similar to that of Theorem 6.1 and is based on Lemmas 5.1 and 6.2.
The only difference is that we do not use conditional blocks to describe
the colorings of all blocks of $\pi$. Thus, the contributions of various $[n]$-colorings 
of $\pi$ to the limit laws reduce to those corresponding to $[r]$-colorings 
and are described in terms of numbers $u(\pi_i,f)$, where $i\in [k]$ and $f\in F_{r}(\pi)$.  
Namely, we have
\begin{eqnarray*}
\lim_{n\rightarrow \infty}v_{0}(\pi,n)&=&
\lim_{n \rightarrow \infty}\left(\frac{1}{n^{k}}
\sum_{f\in F_r(\pi)}n_{f_{1}}u(\pi_1, f)\ldots n_{f_{k}}u(\pi_k,f)\right)\\
&=&
\sum_{f\in F_{r}(\pi)}d_{f_1}u(\pi_1,f)\ldots d_{f_{k}}u(\pi_k,f)\;\;=\;{\rm Tr}(B_{0}(\pi))
\end{eqnarray*}
as $n\rightarrow \infty$, where $\pi\rightarrow B_{0}(\pi)$ is the 
function defined by Definition 5.2, which proves our assertion.
\hfill $\blacksquare$

\section{Matricial semicircle distributions}

The results of Section 6 lead to combinatorial formulas for the asymptotic moments
in the corresponding central limit theorems. In this Section we are going to 
express the limits in terms of their Cauchy transforms represented in the form of
continued fractions. They play the role of the (standard and tracial) `matricial 
semicircle distributions'.

For that purpose
let us recall definitions of certain convolutions of 
distributions, or more generally, of probability measures.
If $F_{\mu}$ is the reciprocal Cauchy transform of some
probability measure $\mu\in {\mathcal M}_{{\mathbb R}}$, then the K-transform 
of $\mu$ is given by $K_{\mu}(z)=z-F_{\mu}(z)$. The boolean 
additive convolution $\mu\uplus \nu$ can be defined by the equation
$$
K_{\mu\uplus \nu}(z)=K_{\mu}(z)+K_{\nu}(z)
$$
where $\mu, \nu \in {\mathcal M}_{{\mathbb R}}$ and $z\in {\mathbb C}^{+}$, 
respectively. In fact, this equation shows that the K-transform is 
the boolean analog of the logarithm of the Fourier transform [19].

We will also need another convolution, which reminds the monotone convolution 
[17], called the orthogonal additive convolution and defined by the equation
$$
K_{\mu \vdash \nu}(z)=K_{\mu}(F_{\nu}(z))
$$
where $\mu, \nu \in {\mathcal M}_{{\mathbb R}}$ and $z\in {\mathbb C}^{+}$.
It was introduced in [12], where we showed that the above formula defines a unique 
probability measure on the real line.
Moreover, if $\mu$ and $\nu$ are compactly supported, 
both $\mu\vdash \nu$ and $\mu \uplus \nu$ are compactly supported.

Using these convolutions, we will now define certain important 
continued fractions which converge
uniformly on the compact subsets of ${\mathbb C}^{+}$ to the K-transforms of some 
probability measures $\mu_{i,j}\in {\mathcal M}_{{\mathbb R}}$.

\begin{Lemma}
For given $B\in M_{r}({\mathbb R})$ with nonnegative entries, continued fractions of the form
$$
K_{i,j}(z)=
\cfrac{b_{i,j}}{z-\sum_{k}\cfrac{b_{k,i}}{z-\sum_{p}\cfrac{b_{p,k}}{z-\ldots}}}
$$
where $i,j \in [r]$, converge uniformly on the compact subsets of ${\mathbb C}^{+}$ to 
the K-transforms of some $\mu_{i,j}\in {\mathcal M}_{{\mathbb R}}$ with compact supports.
\end{Lemma}
{\it Proof.}
Let us define a sequence of functions which approximate the $K_{i,j}$.
Namely, set $K_{i,j}^{(0)}(z)=b_{i,j}/z$ for any $i,j \in [r]$,
which are the K-transforms of probability measures on ${\mathbb R}$ for any $i,j$ 
(Bernoulli measures if $b_{i,j}>0$ and $\delta_{0}$ if $b_{i,j}=0$). In order to use
an inductive argument, let us establish the recurrence
$$
K_{i,j}^{(m)}(z)=\frac{b_{i,j}}{z-\sum_{p}K_{p,i}^{(m-1)}(z)}
$$
for $m\geq 1$. If the $K_{p,i}^{(m-1)}$ are the K-transforms of 
some $\mu_{p,i}^{(m-1)}\in {\mathcal M}_{{\mathbb R}}$ for any $i$ and $p$, 
respectively, then the sums $\sum_{p}K_{p,i}^{(m-1)}$ are the K-transforms of some
$$
\mu_{i}^{(m-1)}=\mu_{1,i}^{(m-1)}\uplus \mu_{2,i}^{(m-1)}\uplus \ldots \uplus \mu_{r,i}^{(m-1)}\;\in 
{\mathcal M}_{{\mathbb R}},
$$
and next, the $K_{i,j}^{(m)}$ are the K-transforms of some
$$
\mu_{i,j}^{(m)}=\kappa_{i,j}\vdash \mu_{i}^{(m-1)}\;\in {\mathcal M}_{{\mathbb R}}
$$
where the $\kappa_{i,j}$ are the Bernoulli measures with K-transforms
$K_{i,j}(z)=b_{i,j}/z$, respectively.
It is easy to see that all these measures are compactly supported.
Moreover, the properties of the orthogonal additive convolution (Corollary 5.3 in [10]) say
that the moments of $\mu_{i,j}^{(n)}$ of orders 
$\leq 2m$ agree with the corresponding moments of
$\mu_{i,j}^{(m)}$ for any $n>m$ and any given $i,j$.
Therefore, we have weak convergence 
$$
{\rm w}-\lim_{m\rightarrow \infty}\mu_{i,j}^{(m)}=\mu_{i,j}
$$
to some $\mu_{i,j}\in {\mathcal M}_{{\mathbb R}}$ for any $i,j$. 
These measures are also compactly supported since $\sup_{i,j}b_{i,j}$ is finite.
In turn, this implies that 
the corresponding Cauchy transforms (and thus K-transforms) 
converge uniformly to the Cauchy transform (K-transforms) of $\mu_{i,j}$ on compact subsets of 
${\mathbb C}^{+}$. This completes the proof.
\hfill $\blacksquare$\\
\indent{\par}
We are ready to state a theorem, which can be viewed as 
the tracial central limit theorem for matricially free random variables.
\begin{Theorem}
Under the assumptions of Theorem 6.1, the
$\psi_{n}$-distributions of $S_{n}$ converge weakly to the distribution 
given by the convex linear combination
$$
\mu=\sum_{j=1}^{r}d_{j}\mu_{j}
$$
where $\mu_{j}=\mu_{1,j}\uplus \mu_{2,j}\uplus \ldots \uplus \mu_{r,j}$ for each $j=1, \ldots , r$
and $\mu_{i,j}$ is the distribution defined by $K_{i,j}$ for any $i,j$.
\end{Theorem}
{\it Proof.}
By Theorem 6.1, we have combinatorial formulas for the moments $M_{m}$ 
of the limit law in the tracial central limit theorem.
The associated distribution extends to a unique compactly supported
probability measure $\mu$ on the real line since its moments are bounded by
the moments of the Wigner semicircle distribution $\sigma_{a}$ with variance 
$a={\rm sup}_{i,j}b_{i,j}$.
Using the mulitplicative formula (5.4) for $B(\pi)$, 
we can formally write the Cauchy transform of $\mu$ in the form
\begin{eqnarray*}
G_{\mu}(z)&=&\sum_{k=0}^{\infty}M_{2k}z^{-2k-1}\\
&=&
\frac{1}{z}+\sum_{k=1}^{\infty}\left(\sum_{\pi\in \mathcal{NC}_{2k}^{2}}{\rm Tr}(B(\pi)D)\right)
z^{-2k-1}\\
&=&
{\rm Tr}\left((z-K(z))^{-1}D\right),
\end{eqnarray*}
which can be called the `trace formula' for $G_{\mu}$, where
\begin{equation}
K(z)=\sum_{k=1}^{\infty}\sum_{\pi\in\mathcal{NCC}_{2k}^{2}}B(\pi)z^{-2k+1}
\end{equation}
is a diagonal-matrix-valued formal power series. 
Moreover, we will show below that each function $K_{j}$ on its diagonal 
is, in fact, the K-transform of some $\mu_{j}\in {\mathcal M}_{{\mathbb R}}$. 
Then, the formal power series given by the trace formula 
is the Cauchy transform of $\mu$ as a convex linear 
combination of Cauchy transforms of probability 
measures.
In fact, using the definition of $B(\pi)$ and (5.3), we obtain
the equation
\begin{equation}
K(z)=\tau((z-K(z))^{-1}B)
\end{equation}
where $K(z)={\rm diag}(K_{1}(z), \ldots , K_{r}(z))$.
By analogy with the scalar-valued case, we can find its solution in the
form of a continued fraction. Namely, observe that each $K_{j}(z)$
has the form of a formal Laurent series
$$
\frac{c_{-1}}{z}+\frac{c_{-3}}{z^{3}}+ \ldots
$$
for some $c_{-1},c_{-3}, \ldots $,
and therefore the above vector equation can be solved by succesive 
approximations. Namely, we set $\mu_{j}$ to be the (compactly supported) probability
measure associated with the K-transform 
$K_{j}(z)=\sum_{i}K_{i,j}(z)$ for each $j \in [r]$, where
the $K_{i,j}$ are given by Lemma 7.1 for $B=DU$.
These K-transforms solve (7.2).
This, together with the trace formula for $G_{\mu}$, gives
$$
G_{\mu}=\sum_{j=1}^{r}d_{j}G_{\mu_j}
$$
where $G_{\mu_j}(z)=1/(z-K_{j}(z))$ is the Cauchy transform of 
$\mu_{j}$ for $j=1, \ldots ,r$. That completes the proof.
\hfill $\blacksquare$

\begin{Remark}
{\rm The Cauchy transform of each $\mu_j$ can be written as 
a continued fraction of the form
$$
G_{\mu_j}(z)=
\cfrac{1}{z- \sum_{i}\cfrac{b_{i,j}}{z-\sum_{k}\cfrac{b_{k,i}}{z-\sum_{p}\cfrac{b_{p,k}}{z-\ldots}}}}
$$
which converges on the compact subsets of ${\mathbb C}^{+}$.}
\end{Remark}

Next, we state a theorem, which plays the role of the standard
central limit theorem for matricially free random variables.

\begin{Theorem}
Under the assumptions of Theorem 6.2, the $\Phi_{n}$-distributions
of $S_{n}$ converge weakly to the distribution 
$$
\mu_{0}=\mu_{1,1}\uplus \mu_{2,2}\uplus \ldots \uplus \mu_{r,r}
$$
where $\mu_{j,j}$ is the distribution defined by $K_{j,j}$ for each $j$.
\end{Theorem}
{\it Proof.}
By Theorem 6.2, we have combinatorial expressions for the 
limit moments $M_{m}$. The proof is similar to that of Theorem 7.1 
and is based on the trace formula for the K-transform of $\mu_0$
$$
K_{\mu_{0}}(z)={\rm Tr}((z-K(z))^{-1}B_{0})
$$
derived from the definition of the function $b_{0}$, which leads to the equation
for the Cauchy transform
$$
G_{\mu_{0}}(z)=
\frac{1}{z-\sum_{j}K_{j,j}(z)},
$$
which completes the proof.
\hfill $\blacksquare$

\begin{Remark}
{\rm The Cauchy transform $G_{\mu_{0}}$ can be written as a continued fraction of the form
$$
G_{\mu_{0}}(z)=
\cfrac{1}{z- \sum_{j}\cfrac{b_{j,j}}{z-\sum_{i}\cfrac{b_{i,j}}{z-\sum_{k}\cfrac{b_{k,i}}{z-\ldots}}}}
$$
which gives a matricial extension of the continued fraction of the Wigner semicircle distribution.}
\end{Remark}

\section{Decompositions in terms of subordinations}
In the one-dimensional case the limit distributions $\mu_{0}$ and $\mu$ 
are related by Proposition 5.1. In particular, if each  
variance matrix $V(n)$ has identical entries equal to one, 
both central limit theorems (standard and tracial) give the 
Wigner semicircle distribution with variance $1$ (of course, the standard case also 
follows from free probability, whereas the tracial case is related to random matrices).

In this Section we will analyze in more detail the limit distributions for the 
two-dimensional case, namely when each variance matrix $V(n)$ consists of four blocks.
They will be expressed in terms of two-dimensional arrays of distributions.
Finding simple analytic formulas for the corresponding four-parameter Cauchy transforms 
and densities does not seem possible in the general case. However, we shall derive 
decomposition formulas for those measures in terms of s-free additive convolutions 
[12], which gives some insight into their structure (see also [18] for recent results on 
the multivariate case). The s-free additive convolution refers to the subordination property for free additive
convolution, discovered by Voiculescu [24] and generalized by Biane [4]. As shown in [12] and [13],
there is a notion of independence, called {\it freeness with subordination}, or simply {\it s-freeness}, 
associated with the s-free additive convolution and its multiplicative counterpart. 

Recall that the s-free additve convolution of $\mu, \nu \in {\mathcal M}_{{\mathbb R}}$ 
is the unique probability measure $\mu \boxright \nu \in {\mathcal M}_{{\mathbb R}}$
defined by the subordination equation
$$
\nu\,\boxplus\, \mu=\nu\vartriangleright (\mu\, \boxright\, \nu),
$$
where $\vartriangleright$ denotes the monotone additive convolution [17].
Equivalently, the above subordination property can be written in terms of 
Cauchy transforms or their reciprocals. 

Using s-free additive convolutions and the boolean convolution, 
we obtain a decomposition of the free additive convolution of the form
$$
\mu \boxplus \nu=(\mu \boxright \nu) \uplus (\nu \boxright \mu),
$$
which allows us to interpret both s-free additive convolutions appearing here as (in general, non-symmetric) halves of 
$\mu\boxplus \nu$. We find it interesting that the limit distributions in the two-dimensional case
will turn out to be deformations of the free additive convolution of semicircle laws 
implemented by this decomposition. In other words, the subordination 
property and the associated convolutions give a natural framework for studying 
matricial generalizations of the semicircle law. 

For simplicity, it will be convenient to use the indices-free notation for the two-dimensional 
matrix of K-transforms:
$$
\left(
\begin{array}{cc}
a(z) & b(z)\\
c(z) & d(z)
\end{array}
\right)
=
\left(
\begin{array}{cc}
K_{1,1}(z) & K_{1,2}(z)\\
K_{2,1}(z) & K_{2,2}(z)
\end{array}
\right)
$$
and 
$$
A=\left(
\begin{array}{cc}
\alpha & \beta \\
\gamma & \delta
\end{array}
\right)=
\left(
\begin{array}{cc}
\sqrt{b_{1,1}} & \sqrt{b_{1,2}} \\
\sqrt{b_{2,1}} & \sqrt{b_{2,2}}
\end{array}
\right)
=\sqrt{B}
$$
where the square root is interpreted entry-wise.

Moreover, we will distinguish two laws by special notations: 
we denote by $\sigma_{\alpha}$ the Wigner semicircle distribution with 
the Cauchy transform
$$
G_{\sigma_{\alpha}}(z)=\frac{z-\sqrt{z^2-4\alpha^2}}{2\alpha^2},
$$
where the branch of $\sqrt{z^2-4\alpha^2}$ is chosen so that 
$\sqrt{z^2-4\alpha^2}>0$ if $z\in {\mathbb R}$ and $z\in (2\alpha, \infty)$, 
and by $\kappa_{\gamma}$ we denote the Bernoulli law with the
Cauchy transform 
$$
G_{\kappa_{\gamma}}(z)=\frac{1}{z-\gamma^2/z},
$$
i.e. $\sigma_{\gamma}=1/2(\delta_{-\gamma}+\delta_{\gamma})$.

Finally, we will also use the {\it boolean compressions} 
of $\mu\in {\mathcal M}_{{\mathbb R}}$, where $t\geq 0$, 
defined by multiplying its K-transform by $t$, namely we define $T_{t}\mu$ to be the 
(unique) probability measure on ${\mathbb R}$, for which 
$$
K_{T_{t}\mu}=tK_{\mu}.
$$ 
These transformations were introduced and studied in [7] and called `$t$-transformations'
of $\mu$. We allow $t=0$, in which case $T_{0}\mu=\delta_{0}$.
In particular, we shall use two-parameter boolean compressions of
semicircle distributions, $\sigma_{\alpha,\beta}=T_{t}\sigma_{\alpha}$ for $t=(\beta/\alpha)^2$,
with
$$
G_{\sigma_{\alpha,\beta}}(z)=\frac{(2\alpha^2-\beta^2)z-\beta^2 \sqrt{z^2-4\alpha^2}}{(2\alpha^2-2\beta^2)z^{2}+2\beta^{4}}
$$
being their Cauchy transforms, where the branch of the square root is 
the same as in the case of $G_{\sigma_{\alpha}}$.
\begin{Theorem}
If $\alpha, \beta, \gamma, \delta\neq 0$, then the diagonal measures $\mu_{j,j}$ 
defined by $K_{j,j}$, where $1\leq j \leq 2$, have the form               
\begin{eqnarray*}
\mu_{1,1}&=&T_{1/t}(\sigma_{\alpha,\beta}\boxright \sigma_{\delta, \gamma})\\
\mu_{2,2}&=&T_{1/s}(\sigma_{\delta, \gamma}\boxright \sigma_{\alpha, \beta})
\end{eqnarray*}
with the non-diagonal measures given by 
$\mu_{1,2}=T_{t}\mu_{1,1}$ and
$\mu_{2,1}=T_{s}\mu_{2,2}$,
where $t=(\beta/\alpha)^2$ and $s=(\gamma/\delta)^2$,
\end{Theorem}
{\it Proof.}
It is easy to see that the following algebraic relations hold:
$$
a(z)=\frac{\alpha^2}{z-a(z)-c(z)}\;\;\;{\rm and}\;\;\;
d(z)=\frac{\delta^2}{z-b(z)-d(z)},
$$ 
$$
b(z)=\frac{\beta^2}{z-a(z)-c(z)}\;\;\;{\rm and}\;\;\;
c(z)=\frac{\gamma^2}{z-b(z)-d(z)}.
$$
Thus, $b(z)=ta(z)$ and $c(z)=sd(z)$, which gives 
$\mu_{1,2}=T_{t}\mu_{1,1}$ and $\mu_{2,1}=T_{s}\mu_{2,2}$. 
In turn, from the equation for $a(z)$, we get
$$
a(z)=\frac{z-c(z)-\sqrt{(z-c(z))^{2}-4\alpha^2}}{2}=K_{\sigma_{\alpha}}(z-c(z))
$$
and thus $\mu_{1,1}=\sigma_{\alpha}\vdash \mu_{2,1}$. In a similar manner we obtain
$\mu_{2,2}=\sigma_{\delta}\vdash \mu_{1,2}$. Therefore, we arrive at the equations
\begin{eqnarray*}
\mu_{1,1}&=&\sigma_{\alpha}\vdash (T_{s}\sigma_{\delta}\vdash T_{t}\mu_{1,1})\\
\mu_{2,2}&=&\sigma_{\delta}\vdash (T_{t}\sigma_{\alpha}\vdash T_{s}\mu_{2,2})
\end{eqnarray*}
since $T_{t}(\mu\vdash \nu)=(T_{t}\mu) \vdash \nu$. 
In order to express the $\mu_{j,j}$ in terms of s-free additive convolutions, we need to use
the properties of the orthogonal convolution. We have shown in [12] that the moment
of order $k$ of $\mu\vdash \nu$ depends on the moments of orders $\leq k$ of $\mu$ and the moments
of orders $\leq k-2$ of $\nu$. This leads to the conclusion that for any compactly supported 
$\mu, \nu \in {\mathcal M}_{{\mathbb R}}$ and the associated sequence of measures $(\mu\vdash_{m} \nu)$, defined recursively
by 
$$
\mu\vdash_{m}\nu =\mu\vdash (\nu \vdash_{m-1} \mu)\;\;\;{\rm with}\;\;\;
\mu\vdash_{1}\nu=\mu\vdash \nu ,
$$
we have weak convergence ${\rm w}-\lim_{m\rightarrow \infty}\mu\vdash_{m}\nu=\mu\boxright \nu$.
If we take $\mu=T_{t}\omega_{\alpha}$ and $\nu=T_{s}\omega_{\delta}$ (these measures
are compactly supported), we get the desired formulas.
\hfill$\blacksquare$\\
\begin{Corollary}
The measures $\mu_{0}, \mu_1, \mu_2$ can be decomposed as
\begin{eqnarray*}
\mu_{0}&=&
T_{1/t}(\sigma_{\alpha,\beta}\boxright \sigma_{\delta, \gamma})
\uplus
T_{1/s}(\sigma_{\delta,\gamma}\boxright \sigma_{\alpha, \beta})
\\
\mu_1&=&T_{1/t}(\sigma_{\alpha, \beta}\boxright \sigma_{\delta, \gamma})
\uplus
(\sigma_{\delta,\gamma}\boxright \sigma_{\alpha, \beta})
\\
\mu_2&=&
T_{1/s}(\sigma_{\delta, \gamma}\boxright \sigma_{\alpha, \beta})
\uplus
(\sigma_{\alpha, \beta}\boxright \sigma_{\delta,\gamma})
\end{eqnarray*}
where the assumptions and notations are the same as in Theorem 8.1.
\end{Corollary}
{\it Proof.} These decompositions follow immediately from
Theorems 7.1, 7.2 and 8.1.
\hfill $\blacksquare$
\begin{Remark}
{\rm The formulas for the diagonal measures $\mu_{j,j}$
in the proof of Theorem 8.1 remind those for 2-periodic continued fractions if 
we take $t=s=1$. The latter are of the same form, except that the semicircle distributions 
are replaced by much simpler Bernoulli laws. 
Nevertheless, if $t=s=1$, the formulas for the $\mu_{j}$
take a simple form
$$
\mu_{j}=\sigma_{\alpha}\boxplus \sigma_{\delta}
$$
where $j=0,1,2$. By Theorem 7.1 and Corollary 8.1, the same formula holds for $\mu$.
Therefore, all measures $\mu_0,\mu_1,\mu_2,\mu$ can be viewed 
as deformations of the free additive convolution of two semicircle distributions, 
implemented by means of boolean compressions.}
\end{Remark}
Let us consider now the situation in which some of
the numbers $\alpha, \beta, \gamma, \delta$ vanish. 
Suitable formulas can be derived algebraically, as
we did in the proof of Theorem 8.1. However, one can 
also obtain the same results by taking weak 
limits in the 
formulas for the measures $\mu_{i,j}$, using the fact that all
measures involved have compact supports.
For that purpose, let us state a few useful facts
about weak limits which will be of interest to us.
Then, we consider eight cases, to which 
the remaining cases are similar (for instance, 
$\alpha=\beta=0$ is similar to $\gamma=\delta=0$). 
\begin{Proposition}
Let $t=\beta^2/\alpha^2$, where $\alpha, \beta>0$, and let 
$\mu\in {\mathcal M}_{{\mathbb R}}$ be compactly supported.
\begin{enumerate}
\item
If $\alpha\rightarrow 0^{+}$, then
\begin{enumerate}
\item
${\rm w}-\lim \,\sigma_{\alpha}=\delta_{0}$,
\item
${\rm w}-\lim\,(\sigma_{\alpha,\beta})=\kappa_{\beta}$,
\item
${\rm w}-\lim\,(T_{1/t}(\sigma_{\alpha,\beta}\boxright \mu))=\delta_{0}$.
\end{enumerate}
\item
If $\beta\rightarrow 0^{+}$, then 
\begin{enumerate}
\item
${\rm w}-\lim\,\sigma_{\alpha,\beta}=\delta_{0}$,
\item
${\rm w}-\lim\,(T_{t}\mu)=\delta_{0}$,
\item
${\rm w}-\lim\,(T_{1/t}(\sigma_{\alpha,\beta}\boxright \mu))=\sigma_{\alpha}\vdash \mu$.
\end{enumerate}
\end{enumerate}
\end{Proposition}
{\it Proof.}
If $\alpha\rightarrow 0^{+}$, then $K_{\sigma_{\alpha}}(z)\rightarrow 0$, which proves (1a).
Here, as well as in the remaining cases, convergence is uniform on compact subsets of ${\mathbb C}^{+}$.
Moreover, $K_{\sigma_{\alpha,\beta}}=\beta^2/(z-\alpha^2K_{\sigma_{\alpha}})\rightarrow \beta^2/z=K_{\kappa_{\beta}}$, which gives (1b). 
In turn, if $\beta\rightarrow 0^{+}$, then $K_{\sigma_{\alpha,\beta}}(z)=\beta^2/(z-K_{\sigma_{\alpha}}(z))\rightarrow 0$, thus
also ${\rm w}-\lim\,\sigma_{\alpha,\beta}=\delta_{0}$, which proves (2a).
If, in addition $t\rightarrow 0^{+}$, then $T_{t}\mu\rightarrow \delta_{0}$
for any $\mu\in {\mathcal M}_{\mathbb R}$, which proves (2b). Finally, 
$$
T_{1/t}(\sigma_{\alpha,\beta}\boxright \mu)=T_{1/t}(T_{t}\sigma_{\alpha}\vdash (\mu \boxright \sigma_{\alpha,\beta}))=\sigma_{\alpha}\vdash (\mu \boxright \sigma_{\alpha,\beta})
$$
and thus the right hand side tends weakly to $\delta_{0}$ as 
$\alpha\rightarrow 0^{+}$, which gives (1c), and tends weakly
to $\sigma_{\alpha}\vdash (\mu \boxright \delta_{0})=
\sigma_{\alpha}\vdash \mu$, which gives (2c). 
This holds for any $\mu\in {\mathcal M}_{\mathbb R}$,
and we also use the right unit property of $\delta$ with respect to the s-free additive convolution, namely
$\mu \boxright \delta_{0}=\mu$. 
\hfill $\blacksquare$

\begin{Corollary}
If some of the entries of the matrix $A$ vanish, we can 
distinguish eight different cases, for which the 
distributions $\mu_{i,j}$ are given by Table 1.
\end{Corollary}
{\it Proof.}
If $a_{i,j}=0$, then $\mu_{i,j}$ is the Dirac delta, which
easily follows from the algebraic equations for the corresponding K-transforms.
An alternative proof can be given by taking weak limits of the formulas 
of Theorem 8.1, as we proceed with the remaining measures.
Thus, if $\delta\rightarrow 0^{+}$, then
\begin{eqnarray*}
\mu_{1,1}&=&{\rm w}-\lim T_{1/t}(\sigma_{\alpha,\beta}\boxright \sigma_{\delta,\gamma})=T_{1/t}(\sigma_{\alpha,\beta}\boxright \kappa_{\gamma})\\
\mu_{1,2}&=&{\rm w}-\lim (\sigma_{\alpha,\beta}\boxright \sigma_{\delta,\gamma})=
\sigma_{\alpha, \beta}\boxright \kappa_{\gamma},\\
\mu_{2,1}&=&{\rm w}-\lim (\sigma_{\delta,\gamma}\boxright \sigma_{\alpha,\beta})=\kappa_{\gamma}\boxright \sigma_{\alpha, \beta},
\end{eqnarray*}
by (1b) of Proposition 8.1, which proves the first case in Table 1. The remaining cases
are proved in a similar manner.
\hfill $\blacksquare$\\

\begin{table}
\caption{Distributions $\mu_{i,j}$ in the case $A$ has zero entries}
\begin{center}
\begin{tabular}{|c|c|c|c|c|c|c|c|}\hline
$a_{1,1}$ & $a_{1,2}$ & $a_{2,1}$ & $a_{2,2}$ & $\mu_{1,1}$ & $\mu_{1,2}$ & $\mu_{2,1}$ & $\mu_{2,2}$ \\ \hline \hline

$\alpha$ 
& 
$\beta$ 
& 
$\gamma$ 
& 
$0$ 
&
$T_{1/t}(\sigma_{\alpha,\beta}\boxright \kappa_{\gamma})$ 
&  
$\sigma_{\alpha,\beta}\boxright \kappa_{\gamma}$
&
$\kappa_{\gamma}\boxright \sigma_{\alpha,\beta} $
& 
$\delta_{0}$
\\

0 
& 
$\beta$ 
& 
$\gamma$ 
& 0
&
$\delta_{0}$
&                              
\;\,\,$\kappa_{\beta}\boxright \kappa_{\gamma}$
&
$\kappa_{\gamma}\boxright \kappa_{\beta}$\;\;\, 
&  
$\delta_{0}$
\\

$\alpha$ 
& 
$0$ 
& 
$\gamma$ 
& 
$\delta$
&
$\sigma_{\alpha}\vdash \sigma_{\delta,\gamma}$
&  
$\delta_{0}$ 
&   
$\sigma_{\delta,\gamma} $
&                      
$\sigma_{\delta}$
\\

0 
& 
$\beta$ 
& 
0 
& 
$\delta$
&
$\delta_{0}$ 
&  
$\kappa_{\beta}$
&
$\delta_{0}$ 
&  
$\sigma_{\delta}\vdash \kappa_{\beta}$
\\

0 
& 
0 
& 
$\gamma$ 
& 
$\delta$
&
$\delta_{0}$             
& 
$\delta_{0}$
&
$\sigma_{\delta,\gamma}$ 
& 
$\sigma_{\delta}$
\\

$\alpha$ 
& 
0 
& 
0 
& 
$\delta$
&
$\sigma_{\alpha}$ 
&
$\delta_{0}$ 
&
$\delta_{0}$      
&   
$\sigma_{\delta}$
\\

$\alpha$ 
& 
0 
& 
0        
& 
0
&
$\sigma_{\alpha}$ 
&   
$\delta_{0}$
&
$\delta_{0}$      
&   
$\delta_{0}$
\\

0 
& 
$\beta $
&
0 
& 
0
&
$\delta_{0}$      
&   
$\kappa_{\beta}$
&
$\delta_{0}$      
&   
$\delta_{0}$
\\
\hline
\end{tabular}
\end{center}
\end{table}
%\indent{\par}
In the case of arbitrary matrix $A$, finding the four-parameter densities of $\mu_{0}$ and $\mu$ is unwieldy.
Below we shall just consider two special cases, in which we can find nice formulas for these measures
for matrices $A$ of arbitrary dimension.  These two cases are of special interest since they are
associated with (aymptotic) freeness and (asymptotic) monotone independence.
\begin{Proposition}
If $A$ is a square $r$-dimensional matrix with identical positive entries $\alpha_{j}$ 
in the $j$-th row, then
$$
\mu_j=\sigma_{\alpha_{1}}\boxplus \sigma_{\alpha_{2}}\boxplus \ldots \boxplus \sigma_{\alpha_{r}}
$$
for each $j\in [r]$, and the $\mu_j$ coincide with $\mu$ and $\mu_{0}$.
\end{Proposition}
{\it Proof.}
Since the columns of $A$ are identical, 
the functions $K_{i,j}$ are the same for all $j$'s. Denote
them $L_{i}=K_{i,j}$, where $i,j\in [r]$. Moreover,
$$
L_{i}(z)=
\frac{b_i}
{z-\sum_{j=1}^{r}L_{j}(z)}
\;\;\;{\rm and}\;\;\;
\sum_{i=1}^{r}L_{i}(z)
=
\frac{\sum_{i=1}^{r}b_{i}}
{z-\sum_{j=1}^{r}L_{j}(z)}
$$
for any $i\in [r]$. 
Therefore, $\sum_{i=1}^{r}L_{i}(z)$ is the K-transform of the measure
$$
\sigma_{\alpha_{1}}\boxplus \sigma_{\alpha_{2}}\boxplus \ldots \boxplus \sigma_{\alpha_{r}}
$$
and since $K_{\mu_{j}}(z)=\sum_{i=1}^{r}K_{i,j}(z)=\sum_{i=1}^{r}L_{i}(z)$, the proof for
$\mu_j$ is completed. It is then easy to see that we get the same result for
$\mu$ and $\mu_{0}$.
\hfill $\blacksquare$
\begin{Proposition}
If $A$ is a lower-triangular $r$-dimensional matrix with identical positive entries $\alpha_{j}$ 
in the $j$-th row below and on the main diagonal, then
$$
\mu_j= \sigma_{\alpha_{j}}\vartriangleright (\sigma_{\alpha_{j+1}}\vartriangleright(\ldots \vartriangleright \sigma_{\alpha_{r}})\ldots )
$$
for each $j\in [r]$. Moreover, $\mu_{0}=\mu_{1}$ and $\mu$ is the convex linear combination of the $\mu_j$ as in Theorem 7.1.
\end{Proposition}
{\it Proof.}
As in the proof of Proposition 8.2, note that the $\mu_{i,j}$ do not depend on $j$
and thus we can set $L_{i}=K_{i,j}$ for any $i\geq j$.
If $i=r$, we have
$$
L_{r}(z)=\frac{b_{r}}{z-L_{r}(z)},
$$
using the continued fraction for the $K_{r,j}$ of Lemma 7.1. Therefore,
$\mu_{r,j}=\sigma_{\alpha_{r}}$ for any $j\leq r$. Next, we have
$$
L_{k}(z)=\frac{b_{k}}{z-\sum_{i=k}^{r}L_{i}(z)}
$$
which leads to
$$
L_{k}(z)=K_{\sigma_{\alpha_{k}}}(z-\sum_{i=k+1}^{r}L_{i}(z))
$$
which gives the orthogonal decomposition of $\mu_{k,j}$, 
$$
\mu_{k,j}=\sigma_{\alpha_{k}}\vdash (\mu_{k+1,j}\uplus \ldots \uplus \mu_{r,j})
$$
for any $j \leq k<r$. Now, we claim that
$$
\mu_{i,j}\uplus \ldots \uplus \mu_{r,j}=
\sigma_{\alpha_{i}}\vartriangleright(\sigma_{\alpha_{i+1}}\vartriangleright (\ldots \vartriangleright \sigma_{\alpha_{r}})\ldots )
$$
for any $1 \leq j \leq i \leq r$. Clearly, it holds for $i=r$ and any $j\leq r$ since we have already shown that $\mu_{r,j}=\sigma_{\alpha_{r}}$ for any $j \leq r$. Suppose now that this formula holds for
$i>k$ and any $j\leq i$. We will show that it holds for $i=k$ and any $j\leq k$. Using the orthogonal decomposition of $\mu_{k,j}$ given above and the inductive assumption, we obtain
$$
\mu_{k,j}=\sigma_{\alpha_{k}}\vdash (\sigma_{\alpha_{k+1}}\vartriangleright(\ldots \vartriangleright \sigma_{\alpha_{r}})\ldots ).
$$
However, for any $\mu,\nu\in {\mathcal M}_{{\mathbb R}}$, we have
a simple relation
$$
(\mu\vdash \nu) \uplus \nu = \mu \vartriangleright \nu
$$
which gives
$$
\mu_{k,j}\uplus \ldots \uplus \mu_{r,j}=\sigma_{\alpha_{k}}\vartriangleright (\sigma_{\alpha_{k+1}}\vartriangleright (\ldots \vartriangleright \sigma_{\alpha_{r}})\ldots )
$$
and the desired expression for $\mu_{j}$.
In a similar manner we obtain $\mu_0$ and $\mu$.
\hfill $\blacksquare$
\begin{Example}
{\rm If $\alpha=\delta\neq 0$ and $\beta=\gamma=0$, then we can use Table 1
to obtain $\mu_{0}=\sigma_{\alpha}\uplus \sigma_{\alpha}$, which is 
the arcsine law with Cauchy transform $G_{\mu_{0}}(z)=1/\sqrt{z^2-\alpha^2}$
whereas $\mu$ is the Wigner semicircle distribution $\sigma_{\alpha}$.
In turn, if $\alpha=\delta=0$ and $\beta=\gamma\neq 0$, then $\mu_{0}=\delta_{0}$,
whereas $\mu=\sigma_{\beta}$.}
\end{Example}
\section{Weighted binary trees and Catalan paths}

In thie Section we show how to express the limit distributions in 
terms of walks on weighted binary trees, or equivalently, in terms of weighted Catalan paths.
The binary tree serves here as an example of the strongly matricially free Fock space. 

The usual framework which gives a description of distributions in 
terms of walks on graphs is the following. Let $W(n)$ denote the set 
of root-to-root walks of lenght $n$ on a rooted graph $({\mathcal G},e)$ 
and let $\mu$ be the spectral distribution of $({\mathcal G},e)$, i.e. 
the distribution given by the moments of the adjacency matrix
$A({\mathcal G})$ in the state $\varphi$ associated with the vector $\delta_{e}$
on the space of square intergrable functions on the set $V({\mathcal G})$
of the vertices of ${\mathcal G}$. Then the $n$-th moment of $A({\mathcal G})$ in the state $\varphi$ 
is equal to the cardinality of the set $W(n)$.
In particular, it is well known that the moments of the Wigner semicircle 
distribution of variance $1$ can be expressed
in terms of walks on the half-line $({\mathbb T}_{1}, e)$ 
with the first vertex denoted by $e$ and chosen as the root. 

For many distributions we have to use a more general framework, in which the moments 
of these distributions are expressed in terms of root-to-root (random or, more generally, weighted) walks on some rooted graph, except that to each walk $w$ on this graph we have to assign a real-valued weight $\xi(w)$. Then we can write
$$
M_{\mu}(n)=\sum_{w\in W(n)}\xi(w)
$$
for any $n\geq 1$, where $\mu$ is the considered distribution. 
In particular, we obtain the moments of $\omega_{\alpha}$ for any $\alpha>0$ by putting $\xi(w)=\alpha^{n}$, where $n=|w|$ is the lenght of $w$. 

In the cases which are of interest to us, the weight function
$\xi$ is first defined on the set of edges $E({\mathcal G})$ of ${\mathcal G}$
and then is extended to $W=\bigcup_{n\geq 1}W(n)$ by multiplicativity. Namely,
if we are given a mapping $\xi:\,E({\mathcal G})\rightarrow {\mathbb R}$, we set
$$
\xi(w)=\xi(E_1)\xi(E_2)\ldots \xi(E_n)
$$
where $w=(E_1,E_2, \ldots , E_n)$ and $E_1,E_2, \ldots , E_n$ are the edges of
$w$ (we choose to describe walks on graphs as sequences of edges).
Such extension, by abuse of notation denoted also by $\xi$,
will be called {\it multiplicative}.

For instance, it is easy to see that the moments of 
$\mu=\sigma_{\alpha} \boxplus \sigma_{\delta}$ can be expressed in this form.
It is enough to take the free product of two half-lines, which is
the binary tree $({\mathbb T}_{2},e)$ with root $e$. 
Let us color this graph in the natural way, namely
each edge which belongs to a copy of the first half-line is colored by $1$
and each edge which belongs to a copy of the second half-line is colored by $2$.
Then the above formula holds for the moments of $\mu$
if we take $\xi(E)=\alpha$ whenever $E$ is colored by $1$ and 
$\xi(E)=\delta$ whenever $E$ is colored by $2$.

\begin{figure}
\unitlength=1mm
\special{em.linewidth 0.5pt}
\linethickness{0.5pt}
\begin{picture}(140.00,60.00)(-25.00,5.00)

\put(10.00,10.00){\line(1,2){5.00}}
\put(20.00,10.00){\line(-1,2){5.00}}
\put(30.00,10.00){\line(1,2){5.00}}
\put(40.00,10.00){\line(-1,2){5.00}}

\put(50.00,10.00){\line(1,2){5.00}}
\put(60.00,10.00){\line(-1,2){5.00}}
\put(70.00,10.00){\line(1,2){5.00}}
\put(80.00,10.00){\line(-1,2){5.00}}

\put(15.00,20.00){\line(2,3){10.00}}
\put(35.00,20.00){\line(-2,3){10.00}}
\put(55.00,20.00){\line(2,3){10.00}}
\put(75.00,20.00){\line(-2,3){10.00}}

\put(25.00,35.00){\line(4,3){20.00}}
\put(65.00,35.00){\line(-4,3){20.00}}

\put(10.00,10.00){\circle*{1.50}}
\put(20.00,10.00){\circle*{1.50}}
\put(30.00,10.00){\circle*{1.50}}
\put(40.00,10.00){\circle*{1.50}}
\put(50.00,10.00){\circle*{1.50}}
\put(60.00,10.00){\circle*{1.50}}
\put(70.00,10.00){\circle*{1.50}}
\put(80.00,10.00){\circle*{1.50}}

\put(10.00,10.00){\line(-1,-2){3.00}}
\put(10.00,10.00){\line(1,-2){3.00}}
\put(20.00,10.00){\line(-1,-2){3.00}}
\put(20.00,10.00){\line(1,-2){3.00}}
\put(30.00,10.00){\line(-1,-2){3.00}}
\put(30.00,10.00){\line(1,-2){3.00}}
\put(40.00,10.00){\line(-1,-2){3.00}}
\put(40.00,10.00){\line(1,-2){3.00}}

\put(50.00,10.00){\line(-1,-2){3.00}}
\put(50.00,10.00){\line(1,-2){3.00}}
\put(60.00,10.00){\line(-1,-2){3.00}}
\put(60.00,10.00){\line(1,-2){3.00}}
\put(70.00,10.00){\line(-1,-2){3.00}}
\put(70.00,10.00){\line(1,-2){3.00}}
\put(80.00,10.00){\line(-1,-2){3.00}}
\put(80.00,10.00){\line(1,-2){3.00}}

\put(15.00,20.00){\circle*{1.50}}
\put(35.00,20.00){\circle*{1.50}}
\put(55.00,20.00){\circle*{1.50}}
\put(75.00,20.00){\circle*{1.50}}

\put(25.00,35.00){\circle*{1.50}}
\put(65.00,35.00){\circle*{1.50}}

\put(45.00,50.00){\circle*{1.50}}

\put(33.00,44.00){\footnotesize $\alpha$}
\put(55.00,44.00){\footnotesize $\delta$}

\put(16.00,27.00){\footnotesize $\alpha$}
\put(31.00,27.00){\footnotesize $\gamma$}
\put(57.00,27.00){\footnotesize $\beta$}
\put(71.00,27.00){\footnotesize $\delta$}

\put(9.00,14.00){\footnotesize $\alpha$}
\put(18.00,14.00){\footnotesize $\gamma$}
\put(28.50,14.00){\footnotesize $\beta$}
\put(39.00,14.00){\footnotesize $\delta$}
\put(49.00,14.00){\footnotesize $\alpha$}
\put(58.00,14.00){\footnotesize $\gamma$}
\put(69.00,14.00){\footnotesize $\beta$}
\put(79.00,14.00){\footnotesize $\delta$}

\put(5.00,6.00){\footnotesize $\alpha$}
\put(12.00,6.00){\footnotesize $\gamma$}
\put(15.50,6.00){\footnotesize $\beta$}
\put(22.00,6.00){\footnotesize $\delta$}
\put(25.50,6.00){\footnotesize $\alpha$}
\put(32.00,6.00){\footnotesize $\gamma$}
\put(36.00,6.00){\footnotesize $\beta$}
\put(42.50,6.00){\footnotesize $\delta$}

\put(45.50,6.00){\footnotesize $\alpha$}
\put(52.00,6.00){\footnotesize $\gamma$}
\put(56.00,6.00){\footnotesize $\beta$}
\put(62.50,6.00){\footnotesize $\delta$}
\put(66.00,6.00){\footnotesize $\alpha$}
\put(72.00,6.00){\footnotesize $\gamma$}
\put(76.00,6.00){\footnotesize $\beta$}
\put(83.00,6.00){\footnotesize $\delta$}

\end{picture}
\caption{Binary tree with a matricial weight function}
\end{figure}
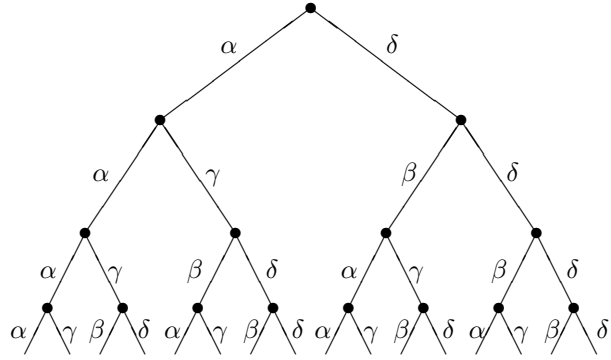
We will demonstrate below that we can express our distributions $\mu_{0},\mu_1,\mu_2$ in a similar form (in particular, we can use the binary tree), except that the weight function $\xi$ 
will depend on all four parameters which appear in the matrix $A$. 
Before formulating the theorem, let us introduce special weight
functions related to matricial freeness.
\begin{Definition}
{\rm Let ${\mathbb T}_{r}$ be a $r$-ary rooted tree with root $e$ and let $A\in M_{r}({\mathbb R})$.
The weight function $\xi:\,E({\mathbb T}_{p})\rightarrow {\mathbb R}$
which assigns the entries of $A$ to the edges of ${\mathbb T}_{r}$ is called {\it matricial}
if, for any pair of edges $E_1,E_2\in E({\mathbb T}_{r})$, incident on the same vertex and
such that $E_1$ is the `father' of $E_2$, the following implication holds: 
$$
\xi(E_1)=a_{i,j}\; {\rm for}\; {\rm some}\; i,j\; \Rightarrow \;\xi(E_2)=a_{k,i}\; 
{\rm for}\; {\rm some}\; k.
$$
The unique multiplicative extension of this weight function to the set of all 
walks on ${\mathbb T}_{r}$ will also be called {\it matricial}.}
\end{Definition}
We specialize to $p=2$ and the binary tree. Note that
any matricial weight function $\xi$ on the binary tree is 
uniquely determined (up to equivalence) by the set 
of those entries of the matrix $A$ which are assigned to the set $\{E_1,E_2\}$ of two edges incident on the root of tree, called the {\it initial weights}. In order to establish a connection with 
our limit distributions, we will assume, as in Section 8, that $A$ is the `square root' of 
$B=DU$ of Section 6. Then, in particular, the binary tree with the 
matricial weight function associated with $A$ and the initial weights $\{\alpha ,\delta\}$, 
shown in Figure 3, describes $\mu_{0}$ as we show below.

Finally, recall after [1,12] that if $({\mathcal G}_{1},e_1)$ and
$({\mathcal G}_{2},e_2)$ are two locally finite simple graphs and $\mu_1$ and $\mu_2$ are
the asscociated spectral distributions, then the s-free product
of ${\mathcal G}_{1}$ and ${\mathcal G}_{2}$ (in that order) 
can be interpreted as this half of the free product 
${\mathcal G}_{1}*{\mathcal G}_{2}$ which `begins' (starting from the root)
with a copy ${\mathcal G}_{1}$. Moreover, the associated 
spectral distribution is given by $\mu_1 \boxright \mu_2$
Let us add that a similar result holds in the multiplicative case:
both the s-free multiplicative convolution and the associated s-free
loop product of graphs were introduced in [13].

Below we shall use the s-free product of half-lines, which are (left and right) halves
of the binary tree.
\begin{Theorem}
Let $\xi_{0}, \xi_1, \xi_2$ be the multiplicative matricial weight functions 
on the set of walks on ${\mathbb T}_{2}$ associated with matrix $A$ and the 
initial weights $\{\alpha,\delta\}$, $\{\beta, \delta\}$ and $\{\alpha, \gamma\}$,
respectively. Then
\begin{equation}
M_{\mu_{j}}(n)=
\sum_{w\in W(n)}
\xi_{j}(w)
\end{equation}
for $j=0,1,2$ and any $n\in {\mathbb N}$, where $W(n)$ denotes the set of root-to-root walks on
${\mathbb T}_{2}$ of lenght $n$.
\end{Theorem}
{\it Proof.}
We need to translate the result of Corollary 8.1 to the language of graphs.
It is well-known that the moments of $\sigma_{\alpha}$ can be interpreted
in terms of walks on the half-line with the weight $\alpha$ assigned to
each edge. The boolean compression $T_{t}$ of $\sigma_{\alpha}$ changes only 
the weight assigned to the edge incident on the root, namely it multiplies it by 
$\sqrt{t}=\beta/\alpha$, which can be illustrated as\\
\unitlength=1mm
\special{em.linewidth 0.5pt}
\linethickness{0.5pt}
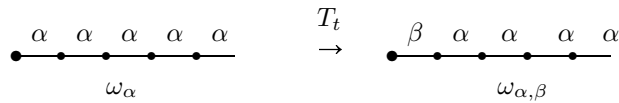
\begin{picture}(140.00,20.00)(-25.00,0.00)

\put(10.00,10.00){\line(1,0){6.00}}
\put(16.00,10.00){\line(1,0){6.00}}
\put(22.00,10.00){\line(1,0){6.00}}
\put(28.00,10.00){\line(1,0){6.00}}
\put(34.00,10.00){\line(1,0){5.00}}

\put(10.00,10.00){\circle*{1.50}}
\put(16.00,10.00){\circle*{1.00}}
\put(22.00,10.00){\circle*{1.00}}
\put(28.00,10.00){\circle*{1.00}}
\put(34.00,10.00){\circle*{1.00}}

\put(12.00,12.00){\footnotesize $\alpha$}
\put(18.00,12.00){\footnotesize $\alpha$}
\put(24.00,12.00){\footnotesize $\alpha$}
\put(30.00,12.00){\footnotesize $\alpha$}
\put(36.00,12.00){\footnotesize $\alpha$}

\put(60.00,10.00){\circle*{1.50}}
\put(66.00,10.00){\circle*{1.00}}
\put(72.00,10.00){\circle*{1.00}}
\put(78.00,10.00){\circle*{1.00}}
\put(84.00,10.00){\circle*{1.00}}

\put(60.00,10.00){\line(1,0){6.00}}
\put(66.00,10.00){\line(1,0){6.00}}
\put(72.00,10.00){\line(1,0){6.00}}
\put(78.00,10.00){\line(1,0){6.00}}
\put(84.00,10.00){\line(1,0){5.00}}

\put(62.00,12.00){\footnotesize $\beta$}
\put(68.00,12.00){\footnotesize $\alpha$}
\put(74.50,12.00){\footnotesize $\alpha$}
\put(82.00,12.00){\footnotesize $\alpha$}
\put(88.00,12.00){\footnotesize $\alpha$}

\put(50.00,14.00){\footnotesize $T_{t}$}
\put(50.00,10.00){\footnotesize $\rightarrow$}
\put(22.00,5.00){\footnotesize $\omega_{\alpha}$}
\put(74.00,5.00){\footnotesize $\omega_{\alpha, \beta}$}

\end{picture}\\
\noindent
Now, the s-free additive convolutions of compressed semicircle distributions which appear
in the decompositions of Corollary 8.1, namely
$$
\sigma_{\alpha, \beta}\boxright \sigma_{\delta,\gamma}\;\;\; {\rm and}\;\;\; 
\sigma_{\delta, \gamma}\boxright \sigma_{\alpha,\beta}
$$ 
are the spectral distributions of the s-free products of 
these half-lines (taken in two different orders). 
These turn out to be the spectral distributions of the two halves of 
the binary tree ${\mathbb T}_{2}$ with the weight function defined by the initial set
$\{\beta, \gamma\}$. In order to obtain $\mu_{0}$, we still need to 
apply $T_{1/t}$ and $T_{1/s}$, respectively, to the left and right halves of the tree, 
which amounts to changing the initial weights from $\beta$ and $\gamma$, respectively, to $\alpha$ and $\delta$.
As a result, we obtain the weight function associated with $A$ and the initial set $\{\alpha,\delta\}$. 
This proves the statement concerning $\mu_{0}$ (the corresponding 
weight function on the binary tree is shown in Figure 3). 
The cases of $\mu_1$ and $\mu_2$ are very similar (at the end, a boolean compression is 
applied to only one half of the binary tree).
\hfill $\blacksquare$\\
\indent{\par}
Another geometric interpretation of Corollary 8.1 can be given in terms of Catalan paths. 
In order to define a Catalan path, we begin with a function $f:[2n]\rightarrow [n]$, such that
$f(0)=f(2n)$ and $|f(i)-f(i-1)|=1$ for any $1\leq i \leq 2n$, and then we 
define a {\it Catalan path} as its unique extension $f:[0,2n]\rightarrow [0,n]$ (by abuse of notation,
denoted by the same symbol) obtained by connecting each $(i-1,f(i-1))$ with $(i,f(i))$ 
with a segment, where $1\leq i \leq 2n$. 
Clearly, each Catalan path consists of segments of two types: 
`rises' $R_1,R_2, \ldots , R_n$, and `falls' $F_1,F_2, \ldots , F_n$.
Moreover, to each `rise' $R_{j}$ there corresponds the closest `fall' 
$F_{\tau(j)}$ lying to the right of $R_{j}$ and on the same (vertical) level.

There is a natural mapping from the set $W(2n)$ of walks of lenght $2n$ on the binary tree
and the set $C(n)$ of Catalan paths of lenght $2n$. In order to rephrase Theorem 8.1 in terms
of Catalan paths, we need to take the sets of {\it weighted Catalan paths},
by which we understand pairs $(f,\xi)$, where $f$ is a Catalan path and 
$\xi$ is a real-valued weight function defined on the set of segments of  
$f$. The multiplicative formula
$$
\xi(f)=\xi(R_{1})\ldots \xi(R_{n})\xi(F_{1})\ldots \xi(F_{n})
$$
assigns the corresponding weight to $f$.
Weighted Catalan paths of special type defined below allow us to
rephrase Theorem 8.1.
\begin{Definition}
{\rm A weighted Catalan path $(f,\xi)$ of lenght $2n$
is called {\it matricial} if $\xi$ assigns entries of 
$A\in M_{p}({\mathbb R})$ to the segments of $f$ in such a way that
the following implications holds:
$$
\xi(R_1)=a_{i,j}\;{\rm for}\; {\rm some}\; i,j\;\Rightarrow\;\xi(R_2)=a_{k,i}
\;{\rm for}\; {\rm some}\; k,
$$
$$
\xi(F_1)=a_{i,j}\;{\rm for}\;{\rm some}\; i,j \;\Rightarrow\; \xi(F_2)=a_{j,k}
\;{\rm for}\;{\rm some}\; k,
$$
for any two consecutive `rises' $R_1,R_2$ and two consecutive `falls' 
$F_1,F_2$, and the same weights are assigned to 
$R_{i}$ and $F_{\tau(i)}$ for each $i \in [n]$.}
\end{Definition}
We will consider below the set of matricially 
weighted Catalan paths of lenght $2n$ associated with the matrix $A$.
In order to rephrase Theorem 9.1, using
weighted Catalan paths, we need to restrict the set of weights
which can be assigned to the first segment of each path (by analogy with trees, 
we call them {\it initial weights}). An example of a weighted Catalan path contributing to
$\mu_{0}$ is given in Figure 4.

\begin{figure}
\unitlength=1mm
\special{em.linewidth 0.5pt}
\linethickness{0.5pt}
\begin{picture}(140.00,30.00)(-25.00,5.00)

\put(10.00,10.00){\line(1,1){5.00}}
\put(15.00,15.00){\line(1,1){5.00}}
\put(20.00,20.00){\line(1,1){5.00}}
\put(25.00,25.00){\line(1,-1){5.00}}
\put(30.00,20.00){\line(1,-1){5.00}}
\put(35.00,15.00){\line(1,1){5.00}}
\put(40.00,20.00){\line(1,-1){5.00}}
\put(45.00,15.00){\line(1,-1){5.00}}
\put(50.00,10.00){\line(1,1){5.00}}
\put(55.00,15.00){\line(1,1){5.00}}
\put(60.00,20.00){\line(1,-1){5.00}}
\put(65.00,15.00){\line(1,1){5.00}}
\put(70.00,20.00){\line(1,-1){5.00}}
\put(75.00,15.00){\line(1,-1){5.00}}

\put(10.00,10.00){\circle*{1.00}}
\put(15.00,15.00){\circle*{1.00}}
\put(20.00,20.00){\circle*{1.00}}
\put(25.00,25.00){\circle*{1.00}}
\put(30.00,20.00){\circle*{1.00}}
\put(35.00,15.00){\circle*{1.00}}
\put(40.00,20.00){\circle*{1.00}}
\put(45.00,15.00){\circle*{1.00}}
\put(50.00,10.00){\circle*{1.00}}
\put(55.00,15.00){\circle*{1.00}}
\put(60.00,20.00){\circle*{1.00}}
\put(65.00,15.00){\circle*{1.00}}
\put(70.00,20.00){\circle*{1.00}}
\put(75.00,15.00){\circle*{1.00}}
\put(80.00,10.00){\circle*{1.00}}

\put(10.50,13.00){\footnotesize $\alpha$}
\put(15.50,18.00){\footnotesize $\gamma$}
\put(20.50,23.00){\footnotesize $\delta$}
\put(27.00,23.00){\footnotesize $\delta$}
\put(32.00,18.00){\footnotesize $\gamma$}
\put(35.50,18.00){\footnotesize $\alpha$}
\put(42.00,18.00){\footnotesize $\alpha$}
\put(47.00,13.00){\footnotesize $\alpha$}
\put(50.50,13.00){\footnotesize $\delta$}
\put(55.50,18.00){\footnotesize $\delta$}
\put(62.00,18.00){\footnotesize $\delta$}
\put(65.50,18.00){\footnotesize $\beta$}
\put(72.00,18.00){\footnotesize $\beta$}
\put(77.00,13.00){\footnotesize $\alpha$}

\end{picture}
\caption{A weighted Catalan path associated with $A$}
\end{figure}
\begin{Corollary}
Let $C_{0}(n)$, $C_{1}(n)$ and $C_{2}(n)$ be the sets of
matricially weighted Catalan paths associated with $A$,
with the initial weights $\{\alpha, \delta\}$, $\{\beta,\delta\}$ and $\{\gamma, \alpha\}$, 
respectively. Then
$$
M_{\mu_{j}}(2n)=
\sum_{(f,\xi)\in C_{j}(n)}
\xi(f)
$$
for $j=0,1,2$ and any $n\in {\mathbb N}$.
\end{Corollary}
{\it Proof.}
This is an easy consequence of Theorem 9.1 since there is a natural bijecton between
each $C_{j}(n)$ and the pair $(W(2n),\xi_{j})$ for any $n$.
\hfill $\blacksquare$


\begin{thebibliography}{99}
\bibitem{[1]}
L. Accardi, R. Lenczewski, R. Sa{\l}apata, Decompositions of the free product of graphs,
{\it Infin. Dimens. Anal. Quantum Probab. Relat. Top.} {\bf 10} (2007), 303-334.
\bibitem{[2]}
M. Anshelevich, Free Meixner states, {\it Commun. Math. Phys.} {\bf 276} (2007), 863-899.
\bibitem{[3]}
D. Avitzour, Free products of $C^{*}$-algebras, {\it Trans. Amer. Math. Soc.}
{\bf 271} (1982), 423-465.
\bibitem{[4]}
Ph. Biane, Processes with free increments, {\it Math. Z.} {\bf 227} (1998), 143-174.
\bibitem{[5]}
M. Bo\.{z}ejko, M. Leinert, R. Speicher, Convolution and limit theorems
for conditionally free random variables, {\it Pacific J.Math.} {\bf 175} (1996), 357-388. 
\bibitem{[6]}
M. Bo\.{z}ejko, R. Speicher,
$\psi$-independent and symmetrized white noises'', in
{\it Quantum Probability and Related Topics VI}, Ed. L.~Accardi,
World Scientific, Singapore, 1991, 170-186.
\bibitem{[7]}
M. Bo\.{z}ejko, J. Wysocza\'{n}ski, Remarks on $t$-transformations of measures and convolutions,
{\it Ann.I.H.Poincar\'{e}- PR} {\bf 37}, 6 (2001), 737-761.
\bibitem{[8]}
T. Cabanal-Duvillard, {\it Probabilit\'{e}s libres et calcul stochastique. Application 
aux grandes matrices al\'{e}eatoires}., Ph.D. Thesis, l'Universit\'{e} Pierre et Marie Curie, 1998.
\bibitem{[9]}
T. Cabanal-Duvillard, V. Ionescu,
Un th\'eor\`eme central limite pour des variables al\'eatoires
non-commutatives, {\it C.\ R.\ Acad.\ Sci.\ Paris},
{\bf 325} (1997), S\'erie I, 1117-1120.
\bibitem{[10]}
W.M.~Ching, Free products of von Neumann algebras, {\it Trans. Amer. Math. Soc.}
{\bf 178} (1993), 147-163.
\bibitem{[11]}
U. Franz, Multiplicative monotone convolutions, in {\it Quantum Probability}, Ed. M. Bo\.{z}ejko {\it et al},
Banach Center Publications, Vol. 73, p. 153-166, 2006.
\bibitem{[12]}
R. Lenczewski, Decompositions of the free additive convolution, 
{\it J. Funct. Anal.} {\bf 246} (2007), 330-365.
\bibitem{[13]}
R. Lenczewski, Operators related to subordination for free multiplicative convolutions,
{\it Indiana Univ. Math. J.}, {\bf 57} (2008), 1055-1103.
\bibitem{[14]}
R. Lenczewski, R. Sa{\l}apata, Discrete interpolation between monotone probability and free
probability, {\it Infin. Dimens. Anal. Quantum Probab. Relat. Top.} {\bf 9} (2006), 77-106.
\bibitem{[15]}
R. Lenczewski, R. Sa{\l}apata, Noncommutative Brownian motions associated with Kesten distributions
and related Poisson processes, {\it Infin. Dimens. Anal. Quantum Probab. Relat. Top.} {\bf 11} (2008), 
to appear.
\bibitem{[16]}
N. Muraki, Monotonic independence, monotonic central limit theorem and monotonic law of small numbers,
{\it Infin. Dimens. Anal. Quantum Probab. Relat. Top.} {\bf 4} (2001), 39-58.
\bibitem{[17]}
N. Muraki, Monotonic convolution and monotonic Levy-Hincin formula, preprint,
2001.
\bibitem{[18]}
A. Nica, Multi-variable subordination distribution for free additive convolution, arXiv:math.OA.0810.2571v1
\bibitem{[19]}
R.~Speicher, R.~Woroudi, Boolean convolution, in {\it Free Probability
Theory}, Ed. D. Voiculescu, 267-279, {\it Fields Inst. Commun.} Vol.12, AMS, 1997.
\bibitem{[20]}
D. Voiculescu, Symmetries of some reduced free product $C^{*}$-algebras, Operator Algebras and Their Connections with Topology and Ergodic Theory, Lecture Notes in Mathematics, Vol. 1132, Springer Verlag, 1985, pp. 556-588.
\bibitem{[21]}
D. Voiculescu, Lectures on free probability theory, {\it Lectures on probability theory and statistics}
(Saint-Flour, 1998), 279-349, Lecture Notes in Math. 1738, Springer, Berlin, 2000.
\bibitem{[22]}
D. Voiculescu, Limit laws for random matrices and free products,
{\it Invent. Math.} {\bf 104}(1991), 201-220.
\bibitem{[23]}
D. Voiculescu , K. Dykema, A. Nica, {\it Free random variables}, CRM Monograph
Series, No.1, A.M.S., Providence, 1992.
\bibitem{[24]}
D. Voiculescu, The analogues of entropy and of Fisher's information measure in free probability theory, I,
{\it Commun. Math. Phys.} {\bf 155} (1993), 71-92.
\bibitem{[25]}
E. Wigner, On the distribution of the roots of certain symmetric matrices, 
{\it Ann. Math.} {\bf 67} (1958), 325-327.
\end{thebibliography}
\end{document}